\documentclass{amsart}

\usepackage[mathcal]{euscript}
\usepackage{mathrsfs}
\usepackage{amsfonts,amssymb,amsbsy,amsmath}
\usepackage{graphicx,color}
\usepackage{epsfig}
\usepackage[all,cmtip]{xy}

\newcommand{\R}{\mathbb{R}}

\newcommand{\C}{\mathbb{C}}
\newcommand{\Z}{\mathbb{Z}}

\newcommand{\ba}{\begin{array}}
\newcommand{\ea}{\end{array}}

\newcommand{\RLF}{\mathcal{RLF}}
\newcommand{\RELF}{\mathcal{RELF}}

\newcommand{\Do}{\mathit{Diff_0}}
\newcommand{\Dco}{\mathit{Diff^c_0}}
\newcommand{\Dcno}{\mathit{Diff_{0}^ {c_{n+1}}}}
\newcommand{\Sg}{\Sigma_g}
\newcommand{\Sa}{\Sigma_1}
\newcommand{\cj}{\mathit{conj}}

\theoremstyle{plain}
\newtheorem{t.}{Theorem}
\newtheorem{l.}[t.]{Lemma}
\newtheorem{p.}[t.]{Proposition}
\newtheorem{c.}[t.]{Corollary}

\theoremstyle{definition}
\newtheorem{d.}{Definition}

\newtheorem{r.}{Remark}

\begin{document}

\title{Invariants of  real Lefschetz fibrations}

\author{Nerm\.{\i}n Salepc\.{\i}}
\address{IUniversit\'e de Strasbourg, DŽpartement mathŽmatiques
7, rue Ren\'e Descartes 67084 Strasbourg Cedex, France}
\email{nsalepci@unistra.fr}

\begin{abstract}
In this note we introduce certain invariants of real Lefschetz fibrations. We call these invariants {\em real Lefschetz chains}. 
We prove that if the fiber genus is greater than 1,  then the real Lefschetz chains are complete invariants of real Lefschetz fibrations with only real critical values. 
If however the fiber genus is 1, real Lefschetz chains are not sufficient to distinguish real Lefschetz fibrations.  
We show that by adding a certain binary decoration to real Lefschetz chains, we get a complete invariant.
\end{abstract}
\maketitle

\section{Introduction}
This note is devoted to a topological study of Lefschetz fibrations equipped with certain $\Z_{2}$ actions compatible with the fiber structure. The action is generated by an involution, which is called a {\em real structure}.  Intuitively, real structures are topological generalizations of the complex conjugation on complex algebraic varieties defined over the reals.  Real Lefschetz fibrations appear, for instance, as blow-ups of pencils of hyperplane sections of complex projective algebraic surfaces defined by real polynomial equations.  Regular fibers of real Lefschetz fibrations are compact oriented smooth genus-$g$ surfaces while singular fibers have a single node. The invariant fibers, called the {\em  real fibers},  inherit a real structure from the real structure of the total space.

 The main results of this article concern topological classification of real Lefschetz fibrations and are exhibited in Section~\ref{gzincirbolum} and Section~\ref{zzincirbolum}  in which we treat the cases of fiber genus $g>1$ and $g=1$, respectively.  In Section~\ref{gzincirbolum}, we introduce {\em real Lefschetz chains} and prove that  if  $g>1$, then real Lefschetz chains are complete invariants of  genus-$g$ real Lefschetz fibrations over a disk that have only real critical values  (Theorem~\ref{gucluzincir}).  The case of $g=1$ (elliptic fibrations) is considered  in Theorem~\ref{zayifzincir}, in  Section~\ref{zzincirbolum}. We show that real elliptic Lefschetz fibrations over a disk that have only real critical values  are determined uniquely by their decorated real Lefschetz chains.  Furthermore, in both cases we study extensions of such fibrations to fibrations over a sphere and obtain complete invariants of real Lefschetz fibrations over a sphere that have only real critical values.

It is possible to give a purely combinatorial shape to decorated real Lefschetz chains. We will discuss such combinatorial objects  (which we call  {\em necklace diagrams}) and their applications in  \cite{N3} (see also \cite{D} for other applications of necklace diagrams).

The present work is organized as follows. In Section~\ref{def}, we settle the definitions and introduce basic notions. 
Section~\ref{erlf} is devoted to the topological classification of  equivariant neighborhoods of real singular fibers.
We show that real Lefschetz fibrations around real singular fibers are determined by the pair consisting of the inherited real structure on one of the nearby regular real fibers and the vanishing cycle which is invariant under the action of the real structure. We call such a pair a {\em real code}. 

In Section~\ref{difuzayi},  we compute the fundamental group of the components of the space of real structures on a genus-$g$ surface. 
These computations are applied in Section~\ref{gucluftoplama}  where we define a {\em strong boundary fiber sum}  (that is, the boundary fiber sum of $\C$-marked real Lefschetz fibrations)  and  show that if the fiber genus is greater than 1, then the strong boundary fiber sum is well-defined. 
As already mentioned above, Section~\ref{gzincirbolum}  gives  a topological classification of genus-$g>1$ Lefschetz fibrations.

Because of the different geometric nature of the surfaces of genus $g>1$ and $g=1$, we apply slightly different techniques to deal with the case of $g=1$.
In Section~\ref{ftoplama}, we define a {\em boundary fiber sum} of non-marked real elliptic Lefschetz fibrations. We observe that the boundary fiber sum is not always well-defined. This observation leads to a decoration of real Lefschetz chains.  In the last section, we introduce {\em decorated real Lefschetz chains} and prove that they  are complete invariants of real elliptic Lefschetz fibrations.  We also study  extensions of such fibrations to fibrations over a sphere.

Let us note that real Lefschetz chains are, indeed, sequences of real codes each of which is associated to a neighborhood of a real singular fiber. Obviously, each real Lefschetz fibration with real critical values defines a real Lefschetz chain which is, by definition, invariant of the fibration.  The natural question to ask is to what extent real Lefschetz chains determine the fibration. This note explores an answer to this question.

\textbf{Acknowledgements.}
This work is extracted from my thesis. I would like to express my gratitude to my supervisors  Sergey Finashin and Viatcheslav Kharlamov  for their guidance throughout  my research.


\section{Basic definitions}\label{def}
Throughout the paper $X$ will stand for a compact connected
oriented smooth 4-manifold and $B$ for a compact connected oriented smooth 2-manifold.

\begin{d.}
A {\em real structure}  $c_X$ on a smooth 4-manifold $X$ is an orientation preserving involution, $c_X^2=id$, such that the set  of fixed
points, $Fix(c_X)$, of $c_X$ is empty or of the middle dimension.

Two real structures $c_X$ and $c'_X$ are considered
{\em equivalent} if there exists an orientation preserving
diffeomorphism $\psi: X \rightarrow X$ such that $\psi \circ c_X=
c'_X \circ \psi$.

A real structure $c_B$ on a smooth 2-manifold $B$ is an orientation reversing
involution $B\to B$. Such structures are similarly considered up
to conjugation by orientation preserving diffeomorphisms of $B$.
\end{d.}

The above definition mimics the properties of the standard complex conjugation on complex
manifolds. Actually, around a fixed point every real structure
defined as above behaves like  complex conjugation.

We will call a manifold together with a real structure a {\em real manifold} and the fixed point set  the {\em real part}.

\begin{r.}\label{reelyapi}
It is well known that for given $g$ there is a finite number of
equivalence classes of {\em real structures} on a genus-$g$ surface 
$\Sigma_g$. These classes can be distinguished by their
{\em types} and the number of real components. 
 Namely, one distinguishes two types of real structures: separating and
non-separating. A real structure is called {\em separating} if the
complement of its real part has two connected components,
otherwise we call it {\em non-separating} (indeed, in the first
case the quotient surface $\Sigma_g/c$ is orientable while in the
second case it is not). The number of real components of a real structure (note
that the real part forms the boundary of $\Sigma_g/c$),  can be at
most $g+1$. This estimate is known as {\em Harnack inequality}
\cite{KRV}. By looking at the possible number of connected
components of the real part, one can see that on $\Sigma_g$ there
are $1+[\frac{g}{2}]$ separating real structures and
$g+1$ non-separating ones. A significant property
of the case of genus-1 surfaces is that the number of real components, which can be 0, 1 or 2, is enough to distinguish the
real structures.
\end{r.}

In this article we stick to the following definition of Lefschetz fibrations.

\begin{d.}\label{LF}
A Lefschetz fibration is a
surjective smooth map
$\pi: X\rightarrow B$
 such that:
\begin{itemize}
\item $\pi(\partial X)=\partial B$ and the restriction $\partial
X\to
\partial B$ of $\pi$ is a submersion;
\item
$\pi$ has only a finite number of
critical points (that is, the points where $d\pi$ is degenerate), all the critical points
belong to $X\setminus\partial X $ and their images
 are distinct points of
$B\setminus\partial B$;
\item around each of the critical points one can choose
orientation-preserving charts $\psi: U \rightarrow \C^2$
and $\phi:V \rightarrow \C$ so that $\phi\circ \pi
\circ{\psi}^{-1}$ is given by
$(z_1,z_2)\rightarrow{z_1}^2+{z_2}^2$.
\end{itemize}
\end{d.}

When we want to specify the genus of the non-singular fibers, we prefer calling them
{\em genus-$g$ Lefschetz fibrations}.
 In particular, we will use the term {\em elliptic Lefschetz fibrations} when the genus is equal to
one. For each integer $g$, we will fix a closed oriented surface of genus $g$, which will serve as a model for the fibers, and denote it by
$\Sigma_g$.  In what follows we will always assume that a Lefschetz
fibration is {\em relatively minimal}; that, is none of its fibers
contains a self intersection -1 sphere. This is not restrictive (if $g\geq1$)
since any self intersection -1 sphere can be blown down while
preserving the projection a Lefschetz fibration.

\begin{d.}\label{RLF}
A real structure on a Lefschetz fibration $\pi: X \rightarrow B$ is a pair of real structures $(c_X, c_B)$ of $X$ and $B$ such that the following diagram commutes

 \begin{displaymath}
\xymatrix{X \ar[r]^{c_X} \ar[d]_{\pi} & X \ar[d]^{\pi} \\
          B \ar[r]^{c_B} & B.}
\end{displaymath}
A Lefschetz fibration equipped with a real structure is called
a {\em real Lefschetz fibration}
and is sometimes referred as $\RLF $.
When the fiber genus is 1, we call it a {\em real elliptic
Lefschetz fibration} (abbreviated $\RELF $).
\end{d.}

\begin{d.}
An {\em $\R$-marked $\RLF $}
is a triple $(\pi, b, \rho)$ consisting of a real Lefschetz fibration 
$\pi:X\rightarrow B$,
a real regular value
$b$ and a diffeomorphism
$\rho: \Sigma_g
\rightarrow
F_b$ such that
$c_X|_{F_{b}}\circ \rho =\rho \circ {c}$
where $c:\Sigma_g\to \Sigma_g$ is a real structure. Let us note that if $\partial B\neq \emptyset$, then $b$ will be chosen in $\partial B$.

A {\em $\C$-marked $\RLF $} is a triple $(\pi, \{b,\bar{b}\},
\{\rho, \bar{\rho}\})$  including a real Lefschetz fibration  $\pi:X\rightarrow B$ a pair of
regular values  $b, \bar{b}=c_{B}(b)$
and a pair of
diffeomorphisms $\rho: \Sigma_g \rightarrow F_b$,
$\bar{\rho}=c _X|_{F_{b}}
\circ \rho:{\Sigma}_g \rightarrow {F}_{\bar{b}}$ where $F_b$ and
${F}_{\bar{b}}=c_X(F_b)$ are the fibers over $b$ and $\bar{b}$,
respectively. As in the case of $\R$-marking, if $\partial B \neq \emptyset$, then we choose $b$ in $\partial B$.
When precision is not needed we will denote $F_b,
F_{\bar{b}}$ by $F$ and $ \bar{F}$, respectively.
\end{d.}

Two real Lefschetz fibrations $\pi: X \rightarrow B$ and $\pi':X' \rightarrow B'$
are said to be {\em isomorphic} if there exist orientation preserving diffeomorphisms
$H: X \rightarrow X'$ and $h: B \rightarrow B'$ such that the following diagram is commutative

 $$
 \xymatrix@!0{
& X \ar@{->}[rr]^H \ar@{->}'[d][dd]_{\pi}
& & X' \ar@{->}^{\pi'}[dd]
\\
X \ar@{->}[ur]^{c_X}\ar@{->}[rr]^H\ar@{->}[dd]_{\pi}
& & X' \ar@{->}[ur]_{c_{X'}}\ar@{->}^{\pi'}[dd]
\\
& B \ar@{->}'[r]^h[rr]
& & B'
\\
B \ar@{->}[rr]^h\ar@{->}[ur]^{c_B}
& & B' \ar@{->}[ur]_{c_{B'}}
.}
 $$

Two $\R$-marked $\RLF s$  are called isomorphic if they are
isomorphic as $\RLF s$ such that
$h(b)=b'$, and the following
diagram is commutative
$$
\xymatrix @!0 @R=2pc @C=3pc {
    F \ar[rr]^{H} \ar[dd]_{c_{X}} && F' \ar[dd]^{c_{X'}} \\
    &\Sigma_g \ar[ru]_{\rho'}\ar[ul]^{\rho} \ar[dd]^c \\
    F \ar[rr]_{H} |!{[ur];[dr]}\hole  && F' \\
    & \Sigma_g \ar[ru]_{\rho'} \ar[lu]^{\rho} .}
    $$

Two $\C$-marked $\RLF s$ are called isomorphic  if they are
isomorphic as $\RLF s$ and
the following diagram is well-defined and
commutative
$$
\xymatrix @!0 @R=2pc @C=3pc {
    F \ar[rr]^{H} \ar[dd]_{c_{X}} && F' \ar[dd]^{c_{X'}} \\
    &\Sigma_g \ar[ru]_{\rho'}\ar[ul]^{\rho} \ar[dd]^{id} \\
    \bar{F} \ar[rr]_{H} |!{[ur];[dr]}\hole  && \bar{F'} \\
    & \bar{\Sigma}_g \ar[ru]_{\bar{\rho'}} \ar[lu]^{\bar{\rho}} .}
    $$

\begin{d.}
A real Lefschetz fibration $\pi:X\to B$ is called {\em directed} if the real part of $(B, c_B)$ is oriented.
(If $c_B$ is separating, then we consider an orientation on the real part inherited from one of the halves $B\setminus Fix(c_B)$.)

Two directed $\RLF s$ are isomorphic if they are isomorphic as $\RLF s$ with the additional condition that the diffeomorphism $h: B\to B$ preserves the chosen orientation on the real part.
\end{d.}

Unless otherwise stated all fibrations we consider are directed.

\begin{r.}
The notion of Lefschetz fibration can be slightly generalized
to cover the case of fibrations whose fibers have non-empty boundary. Then, $X$ turns into a
manifold with corners  and its boundary, $\partial X$, becomes
naturally divided into two parts: the {\em vertical boundary} $\partial^v X$ which is the inverse image $\pi^{-1}(\partial B)$,
and the {\em horizontal boundary}  $\partial^h X$ which is formed
by the boundaries of the fibers. We call such fibrations 
{\em Lefschetz fibrations with boundary}.
\end{r.}


\section{Elementary real Lefschetz fibrations}\label{erlf}

In this section, we classify  real structures on a
neighborhood of a real singular fiber of a real Lefschetz fibration. Such a
neighborhood can be viewed as a Lefschetz fibration over a disc
$D^2$ with a unique critical value $q=0\in D^2$. We call such  a fibration an {\em elementary real Lefschetz fibration}. 
Without loss of generality, we may assume that the real structure on
$D^2$ is the standard one, $\cj$,  induced from $\C\supset D^2$.

Let $\pi: X\to D^2$ be an elementary $\RLF$.
By definition, there exist equivariant local charts
$(U, \phi_U)$, $(V, \phi_V)$ around the critical point
$p\in\pi^{-1}(0)$ and the critical value
$0\in D^2$ respectively
such that $U$ and $V$ are closed discs and $\pi|_U:(U,c_U)\to (V,\cj)$ is equivariantly isomorphic (via $\phi_U$ and
$\phi_V$) to either of $\xi_\pm:(E_\pm,\cj)\to (D_\epsilon,\cj)$, where 
$$E_\pm=\{(z_1, z_2)\in \C^2\, :\, \left|z_1\right|\leq\sqrt{\epsilon},\,\left|z_1^2\pm z_2^2\right|\leq\epsilon^2\}$$ and
$$D_\epsilon=\{t\in \C\, :\, \left|t\right|\leq\epsilon^2 \},\, 0<\epsilon < 1$$ with $\, \xi_\pm(z_1, z_2)=z_1^2\pm z_2^2$.

\begin{figure}[ht]
\begin{center}
   \includegraphics[scale=0.25,trim=0 0 -40 0]{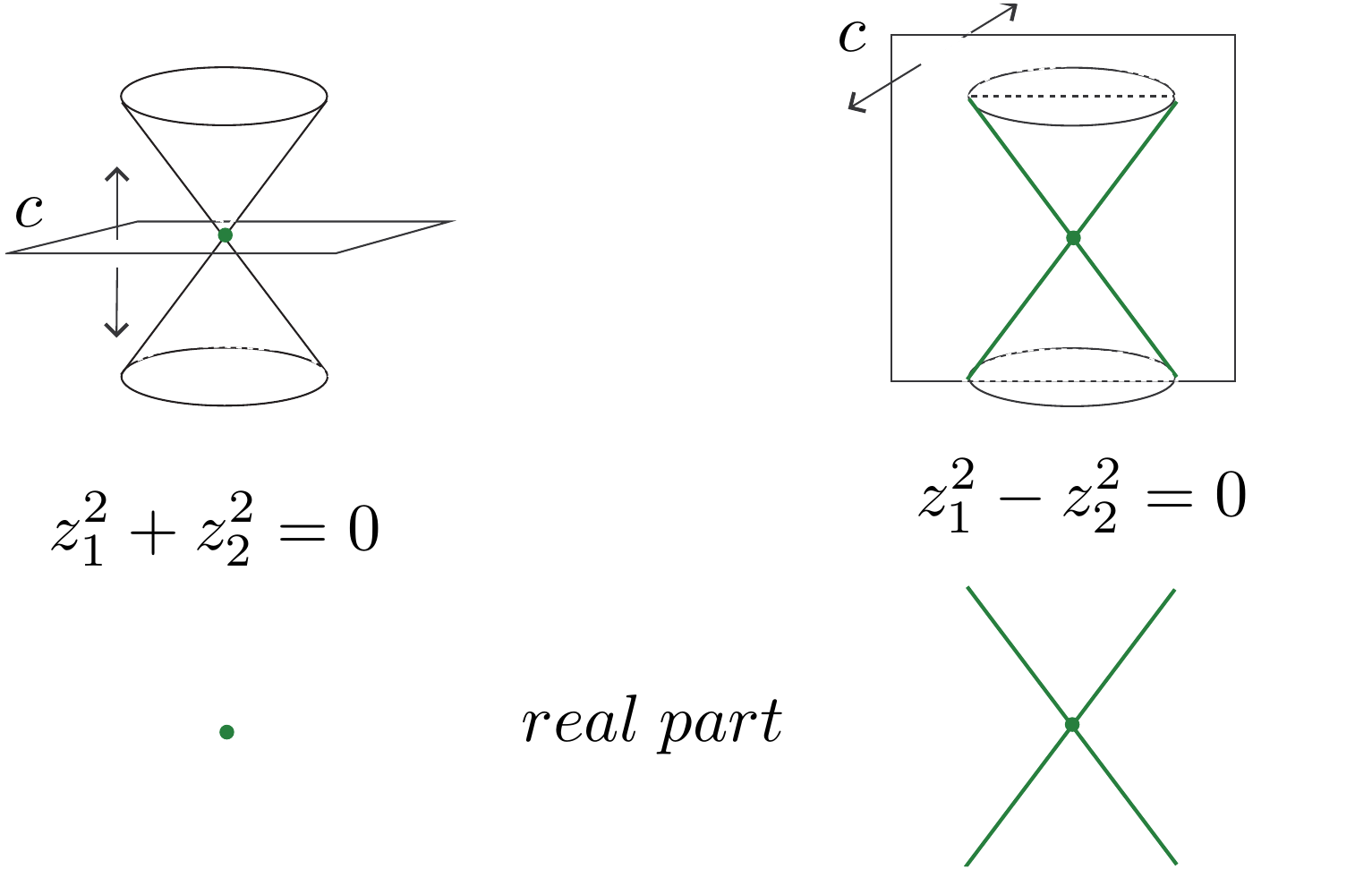}
\caption{\small{Actions of real structures on the singular fibers of $\xi_\pm$.}}
\label{etki}
\end{center}
\end{figure}

The above real local models, $\xi_\pm : E_\pm\to D_\epsilon$, can be seen as two real structures
on the neighborhood of a critical point. These two real structures are not equivalent. The difference can be seen
already at the level of the singular fibers: in the case of $\xi_+$ the two branches are imaginary and they are interchanged by the complex conjugation; in the case of $\xi_-$ the two branches are both real (see Figure~\ref{etki}). 

To understand the action of the real structures on the regular real fibers of $\xi_\pm$, we can use the branched covering defined by
the projection $(z_1, z_2)\to z_1$.
Thus, we have:
\begin{itemize}
\item in the case of $\xi_+$, there are two types of real regular
fibers; the fibers $F_t$ with $t<0$ have no real points, their
vanishing cycles have invariant representatives (that is,
$c(a_t)=a_t$ set-theoretically), and in this case, $c$ acts on the
invariant vanishing cycles as an antipodal involution; the fibers
$F_t$ with $t>0$ has a circle as their real part and this circle
is an invariant (pointwise fixed) representative of the vanishing
cycle; \item in the case of $\xi_-$, all the real regular fibers
are of the same type and the real part of such a fiber consists of
two arcs each having its endpoints on the two different boundary components of the
fiber; the vanishing cycles have invariant representatives,
and $c$ acts on them as a reflection.
\end{itemize}

\begin{figure}[ht]
\begin{center}
   \includegraphics[scale=0.20,trim=0 0 -10 0]{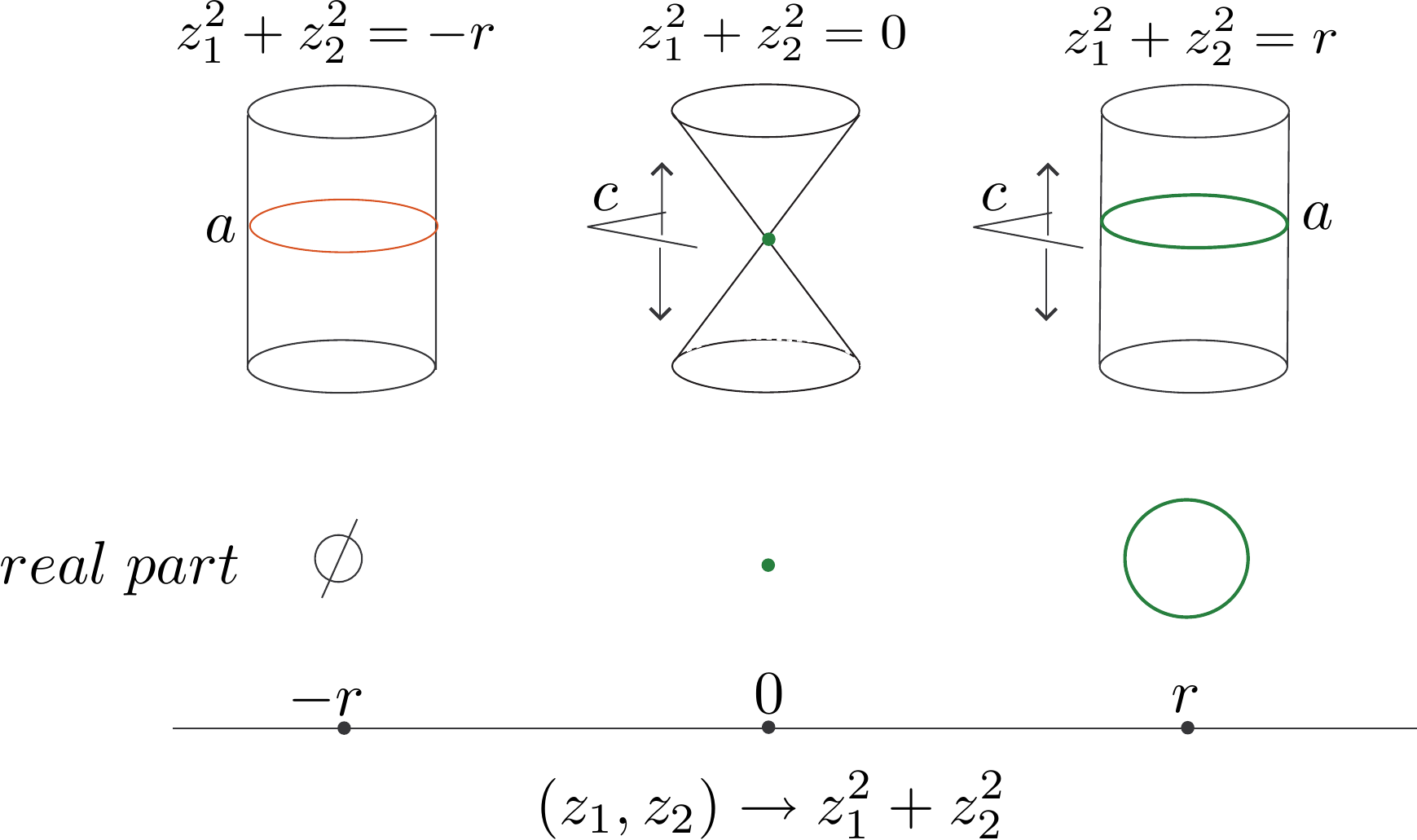}\hspace{.24cm}
   \includegraphics[scale=0.20,trim=0 0 -10 0]{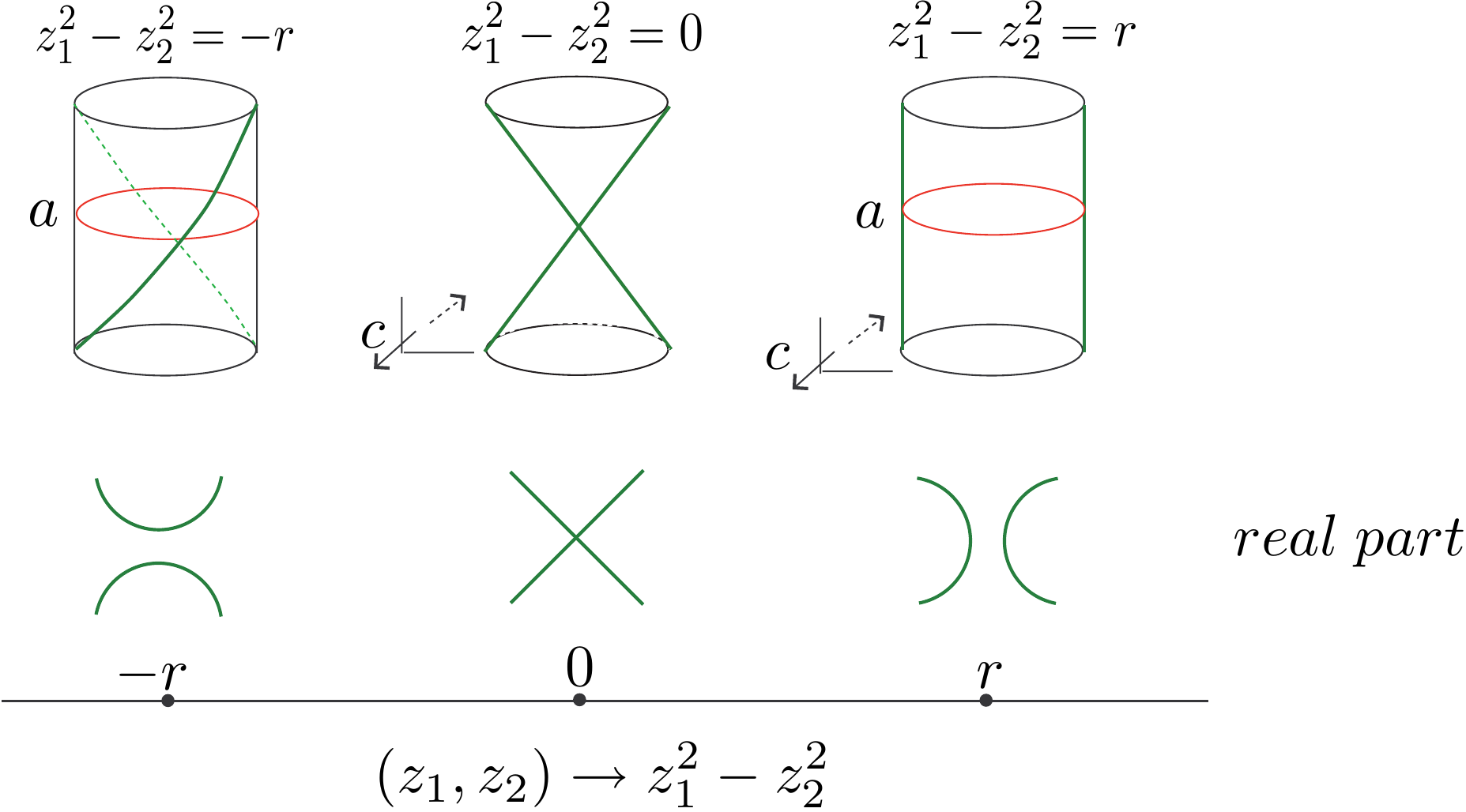}
\caption{\small{Nearby regular fibers of $\xi_\pm$ and the vanishing cycles.}}
\label{etki+}
\end{center}
\end{figure}

Using the ramified covering $(z_1,z_2)\to z_1$, we observe that the horizontal boundary of the
fibration $\xi_\pm$ is equivariantly trivial and has a
distinguished equivariant trivialization. Moreover, since
the complement of $U$ in $\pi^{-1}(V)$ does not contain any
critical point, $X$ can be written as union of two $\RLF s$ with
boundary: one of them, $U\to V$, is isomorphic to $\xi_\pm : E_\pm
\to D_\epsilon$, and the other one is isomorphic to the trivial real fiber bundle $R\to
D_\epsilon$ whose real fibers are equivariantly diffeomorphic to the complement
of an open regular neighborhood of the vanishing cycle $a \subset F_b$. 
The action of the complex conjugation on the
boundary components of the real fibers of  $R\to D_\epsilon$  determines the type of the model, $\xi_\pm: E_\pm
\to D_\epsilon$ glued to  $R\to D_\epsilon$: in the case of
$\xi_+$ it switches the boundary components while in the case of
$\xi_-$ boundary components are preserved (and the complex conjugation acts as a
reflection on each of them).

We use the above decomposition to get  first a classification of directed $\R$-marked elementary $\RLF$,  and then  we discuss the cases of $\C$-marked and non-marked fibrations.

Let $\mathcal{A}^c_g$ denote the set of equivariant isotopy classes of
non-contractible curves on a real surface $(\Sigma_g, c)$ and $\mathcal{V}^c_g$ the set of equivariant isotopy classes of non-contractible
embeddings $\nu:S^1\times I \to \Sigma_g$ such that $c\circ
\nu=\nu$ and  $\mathcal{L}^{\R, c}_g$ the set  of isomorphism classes of directed $\R$-marked elementary real Lefschetz fibrations whose distinguished fiber is identified with $(\Sg, c)$.

We consider the map $\Omega:\mathcal{V}_g^c\to \mathcal{L}_g^{\R, c}$ defined as follows. Let  $[\nu]$ be a class in $ \mathcal{V}_g^c$ with a representative $\nu$.  As $c\circ\nu=\nu$,  the closure,  $\Sigma_g^{\nu}$,  of $\Sigma_g \setminus \nu(S^1\times I)$ inherits a real structure from $(\Sigma_g, c)$.  Let  $R^{\nu}= \Sigma_g^{\nu}\times D_\epsilon\to D_\epsilon$ be the  trivial real fibration with the real structure $c_{R^{\nu}}=(c, \cj):R^{\nu}\to R^{\nu}$ and let  $E^{\nu}_{\pm}\to D_\epsilon$ denote the model $\xi_\pm: E \to D_\epsilon$ whose marked fiber is identified with  $\nu(S^1\times I)$. Depending on the real structure on the horizontal boundary $S^1\times D_\epsilon \to D_\epsilon$ (where the real structure on $S^1\times D_\epsilon$ is taken as $(c_{\partial\Sigma_g^{\nu}}, \cj)$) of $R^{\nu}\to D_\epsilon$, we choose either of $E^{\nu}_{ \pm}\to D_\epsilon$. We then glue $R^{\nu}\to D_\epsilon$ and the suitable model $E^{\nu}_{\pm}\to D_\epsilon$ along their horizontal trivial boundaries to get a fibration in $\mathcal{L}_g^{\R, c}$.

\begin{l.}
$\Omega:\mathcal{V}_g^c\to \mathcal{L}_g^{\R, c}$ is well-defined.
\end{l.}

\noindent {\it Proof:} 
Let $\nu_t :S^1\times I \to \Sigma_g$ be an isotopy between $\nu_{0}$ and $\nu_{1}$.
Then, there exists an equivariant ambient isotopy
$\Psi_t:\Sigma_g\to \Sigma_g$ such that
$\Psi_0=id$ and
$\nu_t=\Psi_t\circ \nu_0$ with $\Psi_t\circ c=c\circ \Psi_t$ for all $t$.
The diffeomorphism,  $\Psi_1$,  induces equivariant diffeomorphisms
$\Psi_1^{R}: R^{\nu_{0}} \to R^{\nu_{1}}$ and $\Psi_1^{E}: E^{\nu_{0}}_{ \pm} \to E^{\nu_{1}}_{ \pm}$ that respect the fibrations and the gluing; thus,
it gives an isomorphism of the images, $\Omega([\nu_{0}])$ and $\Omega([\nu_{1}])$,  as $\R$-marked fibrations. \hfill  $\Box$ \\

Since $c\circ \nu=\nu$, we have $c(\nu(S^1\times \{\frac{1}{2}\}))=\nu(S^1\times \{\frac{1}{2}\}))$. Hence, we can define $\varepsilon:\mathcal{V}_g^c \to \mathcal{A}_g^c$ such that $\varepsilon([\nu])=[\nu(S^1\times\{\frac{1}{2}\})]$. This mapping is two-to-one. Since the monodromy does not depend on the orientation of the vanishing cycle,  there exists a well-defined mapping $\hat{\Omega}$ such that the following diagram commutes

\begin{displaymath}
\xymatrix{ \mathcal{V}_g^c\ar[d]_{\Omega} \ar[r]^{\varepsilon}& \mathcal{A}_g^c \ar[ld]^{\hat{\Omega}}\\
\mathcal{L}_g^{\R, c}. & }
\end{displaymath}

\begin{t.}\label{reelkodR}
$\hat{\Omega}: \mathcal{A}^c_g \to \mathcal{L}^{\R, c}_g$ is a bijection.
\end{t.}

\noindent {\it Proof:}  
As it is discussed in the beginning of the section, any elementary $\RLF$ can be divided equivariantly into  two $\RLF s$ with boundary: an equivariant neighborhood of the critical point (isomorphic to one of the models, $\xi_\pm$), and  the complement of this neighborhood (isomorphic to a trivial real Lefschetz fibration). Such a decomposition defines the equivariant isotopy class of the vanishing cycle. Thus, $\hat{\Omega}$ is surjective.

To show that $\hat{\Omega}$ is injective, let us consider the classes, $[a], [a'] \in \mathcal{A}^c_g$, such that $\hat{\Omega}([a])=\hat{\Omega}([a'])$.  Let $\pi: X\to D_{\epsilon} $ (respectively, $\pi': X'\to D_{\epsilon} $ ) denote the image, $\hat{\Omega}([a])$, of $[a]$ (respectively, the image, $\hat{\Omega}([a'])$, of $[a']$). Since $\hat{\Omega}$ is well-defined, there exist equivariant orientation preserving diffeomorphisms $H: X\to X'$ and $h:D_\epsilon\to D_\epsilon$ such that we have the following commutative diagrams

\begin{center}$
 \xymatrix@!0{
& X\ar@{->}[rr]^H \ar@{->}'[d][dd]_{\pi}
& & X' \ar@{->}^{\pi'}[dd]
\\
X\ar@{->}[ur]^{c_{L^{a}}}\ar@{->}[rr]^H\ar@{->}[dd]_{\pi}
& & X' \ar@{->}[ur]_{c_{L^{a'}}}\ar@{->}^{\pi'}[dd]
\\
& D_\epsilon \ar@{->}'[r]^h[rr]
& & D_\epsilon
\\
D_\epsilon \ar@{->}[rr]^h\ar@{->}[ur]^{\cj}
& & D_\epsilon \ar@{->}[ur]_{\cj}
}
 $
$
\xymatrix @!0 @R=2pc @C=3pc {
    F \ar[rr]^{H} \ar[dd]_{c_{L^{a}}} && F' \ar[dd]^{c_{L^{a'}}} \\
    &\Sigma_g \ar[ru]_{\rho_{\nu'}}\ar[ul]^{\rho_\nu} \ar[dd]^c \\
    F \ar[rr]_{H} |!{[ur];[dr]}\hole  && F' \\
    & \Sigma_g \ar[ru]_{\rho'} \ar[lu]^{\rho} .}
$    \end{center}

Clearly,  $H(\rho(a))$ is equivariantly isotopic to $\rho'(a')$ where $a$ and $a'$ are representatives of $[a]$ and $[a']$, respectively. Moreover, since $H\circ\rho=\rho$, we have $H(\rho(a))=\rho'(a)$, so $\rho'(a)$ is equivariant isotopic to $\rho'(a')$. 

Let $\psi_t:F'\to F'$, $t\in [0,1]$ such that $\psi_0=id$ and $\psi_1(\rho'(a))=\rho'(a')$ and that $\psi_t\circ c'=c'\circ \psi_t$ for all $t\in[0,1]$. Then, $\Psi_t=\rho'^{-1}\circ \psi_t \circ \rho':\Sigma_g\to \Sigma_g$ is the required isotopy between $a$ and $a'$.
\hfill  $\Box$ \\

Theorem~\ref{reelkodR} shows that $c$-equivariant isotopy classes of vanishing cycles  on $(\Sg, c)$ classify  directed $\R$-marked elementary $\RLF s$. To obtain a classification for directed $\C$-marked $\RLF s$, we study the difference between two $\C$-markings. 

\begin{figure}[ht]
\begin{center}
\includegraphics[scale=0.33,trim=0 0 170 0]{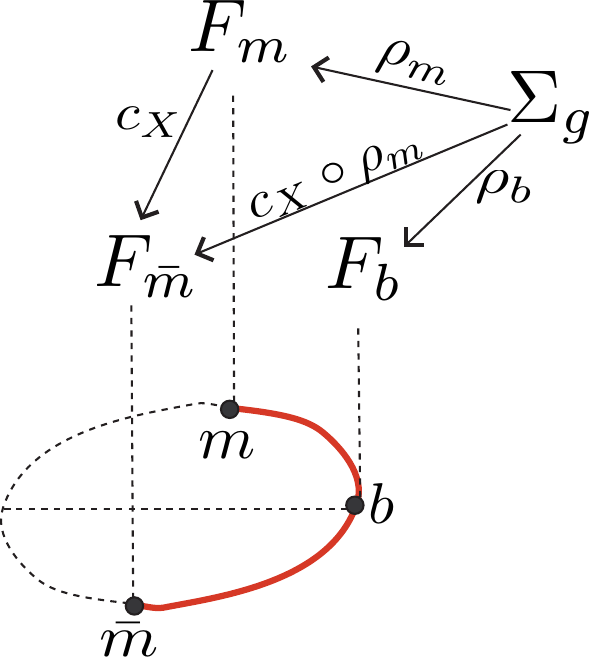}
\caption{Relation between $\R$-marking and $\C$-marking.}
\label{3mbm}
\end{center}
\end{figure}

Let $(\{m, \bar{m}\}, \{\rho_m, c_X\circ \rho_m\})$ be a $\C$-marking on a directed $\RLF$, $\pi:X\to D^2$. The complement, $\partial D^2\setminus \{m,\bar{m}\}$, has two pieces $S_\pm$ (left/ right semicircles) distinguished by the direction. By considering a trivialization of the fibration over the piece of $S_+$ connecting $m$ to the marked real point $b$ (the trivialization over the piece connecting $\bar{m}$ to the real point obtain by the symmetry), we can pull the marking $\rho_m:\Sg\to F_m$ to $F_b$ in order to obtain a marking $\rho_b:\Sg\to F_b$ and a real structure $c=\rho_b^{-1}\circ c_X \circ \rho_b: \Sg \to \Sg$. Any other trivialization results in another marking isotopic to $\rho_b$ and a real structure isotopic to $c:\Sg \to \Sg$.
Hence, a directed elementary $\C$-marked $\RLF$ defines a vanishing cycle defined up to $c$-equivariant isotopy where the real structure $c$ is  also considered up to isotopy.

\begin{d.}
A pair $(c, a)$ of a real structure $c:\Sigma_g\to \Sigma_g$ and a non-contractible simple closed curve $a\in \Sigma_g$ is called a {\em real code}  if $c(a)=a$. 

Two real codes, $(c_0, a_0)$, $(c_1, a_1)$,  are said to be isotopic if there exist a pair of isotopies,  $(c_t, a_t)$,  of real structures and vanishing cycles such that $c_t(a_t)=a_t$,  $\forall  t\in [0, 1]$. Two real codes, $(c_0, a_0)$ and $(c_1, a_1)$, are called conjugate if there is an orientation preserving diffeomorphism $\phi:\Sigma_g\to \Sigma_g$ such that $\phi \circ c_0= c_1 \circ \phi$ and that  $\phi(a_0)$ is isotopic to $a_1$. 

We denote the isotopy class of the real code by $[c, a]$ and the conjugacy class by $\{c, a\}$.  
\end{d.}

\begin{p.}\label{reelkodC}
There is a one-to-one correspondence between the isomorphism classes of directed $\C$-marked elementary $\RLF s$ and the isotopy classes of real codes.
\end{p.}

\noindent {\it Proof:}  Above we discuss how to assign a real code to a directed  $\C$-marked elementary $\RLF$.  It is straightforward to check that this map  is well-defined and surjective.

To show that it is injective,  we consider two isotopy classes $[c_i, a_i]$, $i=1, 2$ such that $[c_1, a_1]=[c_2, a_2]$. Let  $(\pi_1:X_1\to D^2, \{m_1, \bar{m}_1\},\{\rho_{m_1}, \bar{\rho}_{m_1}\})$ and $(\pi_2:X_2\to D^2, \{m_2, \bar{m}_2\},\{\rho_{m_2}, \bar{\rho}_{m_2}\})$ be two directed $\C$-marked elementary $\RLF s$, associated to the classes $[c_1, a_1]$ and $[c_2, a_2]$, respectively. We need to show that $\pi_1$ and $\pi_2$ are isomorphic as directed $\C$-marked $\RLF s$.

Note that we can always choose a representative $c$ for both $[c_1]$ and $[c_2]$ such that $[a_1]=[a_2] \in \mathcal{A}^c_g$. Then,   by Theorem~\ref{reelkodR}, $\pi_1$ is isomorphic to $\pi_2$ as $\R$-marked $\RLF s$. An isomorphism of $\R$-marked $\RLF s$ may not preserve the $\C$-markings; however, it can be modified to preserve them.

\begin{figure}[ht]
\begin{center}
\includegraphics[scale=0.18,trim=0 0 130 0]{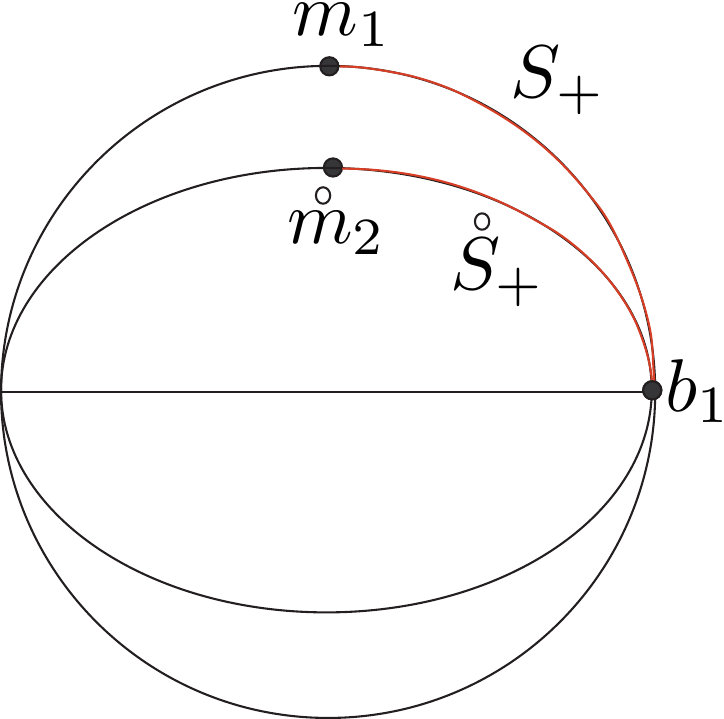}
\caption{The difference of two $\C$-markings.}
\label{3cr}
\end{center}
\end{figure}

Up to homotopy one can identify $X_2$ with a subset, $\stackrel{\circ}{X}_2$, of $X_1$.  Since the difference $X_1\setminus \stackrel{\circ}{X}_2$ has no singular fiber, one can transform the marking $\stackrel{\circ}{m_2}$  of $\stackrel{\circ}{X}_2$ to $m_1$ preserving  the real marking  and the trivializations over the corresponding paths, $S_+$ and $\stackrel{\circ}{S_+}$ (see Figure~\ref{3cr}). This way we get an isomorphism of $\C$-marked $\RLF s$ preserving the isomorphism class of $\R$-marked $\RLF s$.
 
\hfill  $\Box$ \\

For fibrations without marking we allow  $[c, a]$  to change by an equivariant diffeomorphism. Hence, we have the following proposition.

\begin{c.}\label{reelkod}
There is a one-to-one correspondence  between the set of conjugacy classes of  real codes and the set of classes of directed non-marked elementary  real Lefschetz fibrations. \hfill  $\Box$ \\
\end{c.}

\begin{r.}  As the classification of real structures on a genus-$g$ surface is known, it is possible to enumerate the conjugacy classes, $\{c,a\}$ of real codes.  In the case when $a$ is non-separating, there are 6 classes if  $g=1$;  $8g-3$ classes if $g>1$ and is odd;  $8g-4$ classes otherwise. 
The formulas for separating curves are rather bad-looking, so  it is not worth writing them down here. Explicit formulas can be found in \cite{N1}. 
\end{r.}

\begin{r.} \label{fcc} Note that there is no preferable real fiber over the boundary of the disk if the fibration is not directed. Thus, to an elementary non-directed $\RLF$,  we can associate  two real codes,  $(c_{-},a_{-})$,  $(c_{+},a_{+})$, extracted from the ``left" and ``right" real fibers, respectively.  It is a fundamental property of the monodromies of real Lefschetz fibrations that the real structures $c_{-}, c_{+}$ are related 
 by the monodromy such that  $c_{+}\circ c_{-}=t_{a_{-}}=t_{a_{+}}$ (cf. \cite{N2}).
\end{r.}


\section{Equivariant diffeomorphisms and the space of real structures}\label{difuzayi}
In this section we compute the fundamental group of the space of real structures on a genus-$g$ surface.  The computations will  be essential in next sections.

 Let $\mathcal{C}^c(\Sg)$ denote the space of real structures on
$\Sigma_g$ which are isotopic to a fixed real structure $c$, and let $\Do(\Sigma_g)$ denote the group
of orientation preserving diffeomorphisms of $\Sigma_g$ which are
isotopic to the identity. We consider two subgroups of $\Do(\Sg)$: the one,
denoted $\Dco(\Sigma_g)$, consists of those diffeomorphisms
which commute with $c$, and the other, $\Do(\Sigma_g,c)$, 
is the group of diffeomorphisms which are $c$-equivariantly isotopic to the identity. 
The group $\Do(\Sigma_g)$ acts transitively on $\mathcal{C}^c(\Sg)$ by
conjugation. The stabilizer of this action is the group
$\Dco(\Sigma_g)$. Hence, $\mathcal{C}^c(\Sg)$ can be identified
with the homogeneous space $\Do(\Sigma_g)/ \Dco(\Sigma_g)$.

\begin{l.}\label{baglantili}
The space  $\Dco(\Sigma_g)$ is connected for all $c:\Sg\to \Sg$  if  $g>1$, and  for  $c:\Sg\to \Sg$ which has one real component if $g=1$.
\end{l.}

\noindent {\it Proof:}  (We will use different techniques for the cases $g>1$ and $g=1$.)

The case of $g>1$:  we consider the fiber bundle description of conformal structures on $\Sg$, introduced in \cite{EE}.  Let $\mathit{Conf}_{\Sg}$ denote the space of conformal structures on $\Sg$ equipped with $C^\infty$-topology. The group $\Do(\Sigma_g)$ acts on $\mathit{Conf}_{\Sg}$ by composition from right. This action is proper, continuous, and effective; hence, $\mathit{Conf}_{\Sg}\to \mathit{Conf}_{\Sg}/ \Do(\Sg)$ is a principle $\Do(\Sg)$-fiber bundle (cf. \cite{EE}).  The quotient is the Teichmuller space of $\Sg$, denoted $\mathit{Teich}_{\Sg}$.
Note that conformal structures can be seen as equivalence classes of Riemannian metrics with respect to the relation that two Riemannian metrics are equivalent if  they differ by a positive function on $\Sg$. Let $\mathit{Riem}_{\Sg}$ denote the space of Riemannian metrics on $\Sg$. Then, we have the following fibrations

\begin{displaymath}
\xymatrix{ \{u:\Sg\to \R: u>0\}\ar[r] &\mathit{Riem}_{\Sg}  \ar[d]_{p_2}  \\
          \Do(\Sg) \ar[r] & \mathit{Conf}_{\Sg} \ar[d]_{p_1} \\
           &\mathit{Teich}_{\Sg}.}
\end{displaymath}

The real structure $c$ acts on $\Do(\Sg)$ by conjugation. This action can be extended to $\mathit{Conf}_{\Sg}$ and $\mathit{Riem}_{\Sg}$ as follows.
We fix a section $s:\mathit{Teich}_{\Sg}\to \mathit{Conf}_{\Sg}$ of the bundle $p_1$ and we consider a family of diffeomorphisms $\phi^s_{\zeta}:\Do(\Sg)\to p_1^{-1}(\zeta)$ parametrized by $\mathit{Teich}_{\Sg}$ such that $\phi^s_{\zeta}(id)=s(\zeta)$. Let $s(\zeta)=[\mu_x]$ for some Riemannian metric $\mu_x$ on $\Sg$. Then, we define $\phi^s_{\zeta}(f(x))=[\mu_{f(x)}]$ for all $f\in \Do(\Sg)$. 
The action of the real structure, thus, can be written as $c.[\mu_{f(x)}]=[\mu_{c\circ f\circ c(x)}]$. Clearly the definition does not depend on the choice of the representative of the class $[\mu_{f(x)}]$, so the action extends to $\mathit{Riem}_{\Sg}$.

Let $Fix_{\mathit{Conf}_{\Sg}}(c)$ denote the set of fixed points of the action of $c$ on $\mathit{Conf}_{\Sg}$ and $Fix_{\mathit{Riem}_{\Sg}}(c)$ be the set of fixed points on $\mathit{Riem}_{\Sg}$.  
Note that $s(\zeta)=\phi^s_{\zeta}(id) \in Fix_{\mathit{Conf}_{\Sg}}(c)$ for all $ \zeta \in \mathit{Teich}_{\Sg}$. Indeed,  each $[\mu_{f(x)}]$ where $f\in \Dco(\Sa)$ is in $Fix_{\mathit{Conf}_{\Sg}}(c)$. 

The space  $Fix_{\mathit{Conf}_{\Sg}}(c)$ is connected. If  $Fix_{\mathit{Conf}_{\Sg}}(c)$ were disconnected, then the inverse image $Fix_{\mathit{Riem}_{\Sg}}(c)$ would also be disconnected in $\mathit{Riem}_{\Sg}$. However, it is known that $\mathit{Riem}_{\Sg}$ is convex; thus,  $Fix_{\mathit{Riem}_{\Sg}}(c)$ is convex, so  it is connected.  Therefore, $Fix_{\mathit{Conf}_{\Sg}}(c) \cap \Do(\Sg)=\Dco(\Sg)$ is connected since $Fix_{\mathit{Conf}_{\Sg}}(c)$ is a union of sections.

The case of  $g=1$: if $c$ has one real component, then the quotient $\Sa/ c$ is the M\"obius band. The space of diffeomorphisms of the M\"obius band has two connected components  \cite{Hm2}: the identity component  and the component of the diffeomorphism induced (if the M\"obius band is obtained from $I\times I$,  by identifying the points $t\times 0$ with the points  $1-t \times 1, t\in I=[0,1]$)  from the reflection of $I\times I$ with respect to $I\times \frac{1}{2}$. This diffeomorphism is not isotopic to the identity  because before identifying the ends it reverses the orientation of $I\times I$, and  it lifts to a diffeomorphism of $\Sa$ (considered as the obvious quotient of $[-1,1]\times [-1,1]$)  induced from the the central symmetry of $[-1,1]\times [-1,1]$.  This diffeomorphism is not isotopic to the identity on $\Sa$ since it reverses the orientation of the real curve. 

Therefore, we have $$\{f:\Sa/ c\to \Sa/ c:  \textrm{ $\hat{f}: \Sa\to \Sa$ is isotopic to $id$} \}=\{f:\Sa/c\to \Sa/ c: f \cong id\}.$$ The former is identified by $\Dco(\Sa)$ and the latter is connected.
\hfill  $\Box$ \\

\begin{l.}\label{pi1}
For any real structure $c:\Sigma_g\to \Sigma_g$, 

\begin{displaymath}
\pi_1(\Do(\Sigma_g)/ \Do(\Sigma_g, c), id)=
\left\{\begin{array}{ll}
0 &\textrm{if $g>1$ } \\
\Z & \textrm{if $g=1$.}
\end{array} \right.
\end{displaymath}
\end{l.}

\noindent {\it Proof:} 
Note that the subgroup $\Do(\Sg, c)$ acts  on $\Do(\Sg)$  by composition from left. Such an action is free, so $\Do(\Sg) \rightarrow \Do(\Sg) / \Do(\Sg, c)$ is a $\Do(\Sg, c)$-fiber bundle.
The fibers, $\Do(\Sg, c)$, can be identified with the group $\Do(\Sg/c)$ because the lifting of diffeomorphisms of $\Sg/c$ can always be assured by means of the orientation double cover of $\Sg/c$. (Note that  if $c$ is non-separating, then $\Sg/c$ is non-orientable. In this case, $\Do(\Sg/c)$ denotes the space of all diffeomorphisms of $\Sg/c$ and $\Do(\Sg/c)$ is component of the identity.)  
 
Now, we consider the long exact homotopy sequence of this fibration.

$$\begin{array}{c}\cdots \rightarrow \pi_2(\Do(\Sg)) \rightarrow \pi_2(\Do(\Sg) / \Do(\Sg, c))\rightarrow \pi_1(\Do(\Sg, c)) \rightarrow \\
\pi_1(\Do(\Sg)) \rightarrow \pi_1(\Do(\Sg) / \Do(\Sg, c))\rightarrow \pi_0(\Do(\Sg))\rightarrow \cdots \end{array}$$

The case of  $g>1$:  the space $\Do(\Sigma_g)$ is contractible for $g>1$ \cite{EE},  so is $\Do(\Sigma_g/ c)$  \cite{ES}.  Therefore, from the homotopy long exact sequence of the fibration we obtain $\pi_1(\Do(\Sigma_g)/\Do(\Sigma_g,c), id)=0$.

The case of  $g=1$:  it is known that $\Sa$ is deformation retract of $\Do(\Sa)$ \cite{I}, so the space $\Do(\Sa)$ can be considered as a group generated by the rotations which lift to the standard translations on the universal cover.  

To understand  $\Do(\Sigma_g,c)$, we first consider the case when $c$ has two real components. Note that, in this case, the quotient $\Sa/c$  is topologically an annulus, so $\pi_1(\Do(\Sa/ c), id)=\Z$ \cite{I}.
We fix an identification of $\varrho: \C/\Z^2\to \Sa$ such that the real structure $c$ is the one induced from the standard complex conjugation on $\C$. 
We consider the following family of diffeomorphisms 
$$\begin{array}{llllllll}
{R'}_{t}^1:& \C/\Z^2&\to&  \C/\Z^2 & {R'}_{t}^2:& \C/\Z^2&\to& \C/\Z^2\\
&(x+iy)_{\Z^2}&\to &(x+t+iy)_{\Z^2}& &(x+iy)_{\Z^2}&\to &(x+i(y+ t))_{\Z^2}
\end{array}$$
where $t\in [0,1]$ and $(x+iy)_{\Z^2}$ denotes the equivalence class of $x+iy$ in $\C/\Z^2$. 
Clearly ${R'}_{0}^j={R'}_{1}^j=id$ and for each $t\in [0,1]$, ${R'}_{t}^j,  j=1,2$ is isotopic to identity. 
The homotopy 
classes of $R_{t}^1=\varrho\circ {R'}_{t}^{1} \circ \varrho^{-1}$ and $R_{t}^2=\varrho\circ {R'}_{t}^{2} \circ \varrho^{-1}$ form a basis of $\pi_{1}(\Do(\Sa), id)$. Moreover,  with respect to the identification $\varrho$,  each diffeomorphism $R_{t}^1$ is in  $\Do(\Sa,c)$, so the loop $R_{t}^1$ is a generator of  $\pi_1(\Do(\Sa,c),id)$. Thus, from the homotopy exact sequence we get $\pi_1(\Do(\Sa) / \Do(\Sa, c), id)=\Z$.

 If $c$ has no real component,  then the quotient $\Sa / c$ is a Klein bottle, so the group $\Do(\Sa / c)$ is isomorphic to $ S^1$ and is generated by the
rotation which lifts to a translation in the universal cover of the Klein bottle \cite{Hm2}.
Let us now fix an identification $\varrho:\R^2/\Z^2\to \Sa$ such that the real structure $c$ is induced from the real structure 
$$\begin{array}{lcl}
 \R^2/ \Z^2 &\to& \R/\Z^2\\
(x, y)_{\Z^2}&\to&(x+\frac{1}{2},-y)_{\Z^2}.
\end{array}$$ 
The classes of family of diffeomorphisms  $R_{t}^j=\varrho\circ {R'}_{t}^j \circ \varrho^{-1}, \, j=1,2$  where  
$$\begin{array}{llllllll}
{R'}_{t}^1:&\R^2/\Z^2&\to& \R^2/\Z^2 & {R'}_{t}^2:&\R^2/ \Z^2&\to& \R^2/ \Z^2\\
&(x, y)_{\Z^2}&\to &(x+t, y)_{\Z^2}& &(x, y)_{\Z^2}&\to &(x, y+ t)_{\Z^2}.
\end{array}$$  form a basis of $\pi_{1}(\Do(\Sa), id)$.  
Moreover,   with respect to the identification $\varrho$ each diffeomorphism $R_{t}^1$ is in $\Do(\Sa, c)$, and so $R_{t}^1$ is a generator of $\pi_1(\Do(\Sa, c), id)=\Z$. Therefore, we  get $\pi_1(\Do(\Sa) / \Do(\Sa,c), id)=\Z$.

If $c$ has a unique real component,  $C$, then the restriction $f|_{C}$ of $f\in \Do(\Sa,c)$  defines a diffeomorphism of $C$.
Such a restriction defines a fibration, $\Do(\Sa,c) \to \Do(C)$, whose fibers isomorphic to $\Do(\Sa, C)=\{f \in \Do(\Sa,c): f|_C=id\}.$
Note that $\Do(\Sa, C)\cong \Do(\overline{\Sa\setminus C}, \partial)$ where $\overline{\Sa\setminus C}$ denotes the closure of $\Sa\setminus C$ and $\Do(\overline{\Sa\setminus C}, \partial)$ the group diffeomorphisms of $\overline{\Sa\setminus C}$ which are identity on the boundary. 

Topologically $\Sa\setminus C$ is an annulus, so $\Do(\overline{\Sa\setminus C}, \partial)$ is contractible \cite{I}. 
 From the homotopy long exact sequence of the following fibration
$$\xymatrix{\Do(\Sa,\textsl{C}) \ar[r] & \Do(\Sa,c) \ar[d] \\
          & \Do(\textsl{C})}$$
   we get $\pi_k(\Do(\Sa, c), id)\cong \pi_k(\Do(C), id)$, $\forall k$.

Let us now choose an identification $\varrho: \C/\Lambda \to \Sa$ where $\Lambda$ is the lattice generated by $v_1=(\frac{1}{\sqrt2},\frac{1}{\sqrt2})$ and $v_2=(\frac{1}{\sqrt2},-\frac{1}{\sqrt2})$. 
Then, the real structure $c$ can be taken as the one induced from the complex conjugation on $\C$.

We consider $R_i'(t):\C/\Lambda\to \C/ \Lambda$, $t\in [0,1]$
such that
$$\begin{array}{llllllll}
{R'}_{t}^1:&\C/\Lambda&\to& \C/\Lambda & {R'}_{t}^2:&\C/ \Lambda&\to& \C/ \Lambda\\
&(x+iy)_{\Lambda}&\to &(x+t+iy)_{\Lambda} & &(x+iy)_{\Lambda}&\to &(x+i(y+ t))_{\Lambda}.
\end{array}$$

Again,  the classes of $ {R_{t}^j}=\varrho\circ  {R'}_{t}^j \circ \varrho^{-1}, j=1,2$ form a basis for $\Do(\Sa)$ while $ {R_{t}^1}$ can be taken as a generator for $\pi_1(\Do(\Sa, c), id)=\Z$. Therefore,  $\pi_1(\Do(\Sa)/ \Do(\Sa, c), id)=\Z$.

\hfill  $\Box$ \\

\begin{p.}\label{fundgrup}
For any real structure $c:\Sigma_g\to \Sigma_g$,
\begin{displaymath}
\pi_{1}(\mathcal{C}^c(\Sg))=\pi_1(\Do(\Sigma_g)/ \Dco(\Sigma_g), id)=
\left\{\begin{array}{ll}
0 &\textrm{if $g>1$ } \\
\Z & \textrm{if $g=1$.}
\end{array} \right.
\end{displaymath}

\end{p.}

\noindent {\it Proof:} By Lemma~\ref{baglantili},
 $\Dco(\Sg)$ is connected for all real $c:\Sg\to \Sg$, $g>1$ and for the real structure
$c:\Sa\to \Sa$ which has one real component. Hence, in these cases
$\Dco(\Sa)=\Do(\Sa,c)$, so the result follows from Lemma~\ref{pi1}.

In the case when $c:\Sa \to \Sa$ has 2 real components, the space $\Dco(\Sa)$ has two connected components. Note that the diffeomorphism ${R_{\frac12}^2}$ ( induced from the translation, $(x+iy)_{\Z^2}\to (x+i(y+ {\frac{1}{2}}))_{\Z^2}$, on $\C/\Z^2$) is equivariant; however, it is not equivariantly isotopic to the identity. Hence,  $\Dco(\Sa)$ has two components:  the component, $\Do(\Sa, c)$, of the identity and the  component of the rotation ${R_{\frac12}^2}$. (In what follows, we denote ${R_{\frac12}^2}$  by $R_{\frac{1}{2}}$.)

We can identify rotations in $\Do(\Sa)\setminus \Do(\Sa,c) $ with  $S^1$ by letting ${R_{t}^2}\to 2\pi t$. 
Then, rotations in the quotient  $\Do(\Sa)/ \Dco(\Sa)$ are identified with $S^1/_{\theta \sim (\theta+\pi)}$,  so  we have  $\pi_1(\Do(\Sa)/\Dco(\Sa), id)=\Z$.

The case when $c:\Sa\to \Sa$ has no real component can be treated similarly using the identification $\varrho:\R^2/ \Z^2\to \Sa$.
\hfill  $\Box$ \\


\section{Boundary fiber sum of $\C$-marked real Lefschetz fibrations}\label{gucluftoplama}
Let $(D^2, \cj)$ be a real disk with oriented real part. We denote by $S^\pm$ the upper/ lower semicircles of $\partial D^2$. We consider also left/ right semicircles, denoted by $S_{\pm}$, and the quarter-circles $S^\pm_\pm=S^\pm \cap S_\pm$. (Here directions right/ left and  up/ down are  determined by the orientations $D^2$ and of the real part.)
Let $r_\pm$ be the real points of $S_\pm$, and $c_\pm$ the real structures on $F_\pm=\pi^{-1}(r_\pm)$.

\begin{d.} Let $(\pi':X'\to D^2, \{b',\bar{b}'\}, \{\rho', \bar{\rho}'\})$ and  $(\pi:X\to D^2, \{b,\bar{b}\}, \{\rho, \bar{\rho}\})$
be two directed $\C$-marked real Lefschetz fibrations such that
the real structures $c'_{+}$ on $F'_+$ and $c_{-}$ on $F_-$ induce (via the markings) isotopic real
structures on $\Sigma_g$. Then, we define the {\em strong boundary fiber sum} (the boundary fiber
sum of $\C$-marked $\RLF s$)  as follows.

\begin{figure}[ht]
   \begin{center}
   \includegraphics[scale=0.36,trim=0 0 80 0]{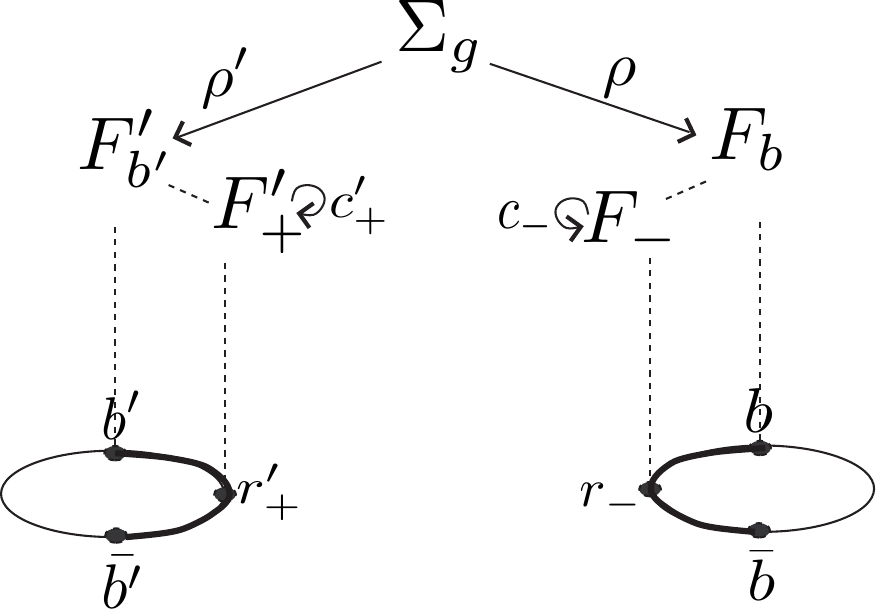}
   \end{center}
 \end{figure}

We choose trivializations of $\pi'^{-1}(S^+_+)$ and $\pi^{-1}(S_-^+)$ such that the pull backs of $c'_{+}$ and $c_{-}$ give the same real structure $c$ on $\Sigma_g$. The trivialization of $\pi'^{-1}(S_+)$ can be obtained as a union  $\Sigma_g \times S^+_{+} \cup \Sigma_g \times S_+^{-}\diagup_{(x,1_+)\sim (c(x),1_-)}$ and similarly  $\pi^{-1}(S_-)=\Sigma_g \times S_-^{+} \cup \Sigma_g \times S_-^{-}\diagup_{(x,-1_+)\sim (c(x),-1_-)}$. The strong boundary fiber sum $X'\natural_{\Sigma_g} X\to D^2\natural D^2$ is, thus, obtained by gluing $\pi'^{-1}(S_+)$ to $\pi^{-1}(S_-)$ via the identity map.
\end{d.}

\begin{r.} (1) In fact, the construction described above  creates a manifold with corners, but there is a canonical way to smooth the corners; hence, the strong boundary fiber sum is the manifold obtained by smoothing the corners.

(2) By definition, the strong boundary fiber sum is associative but not commutative.

(3) The strong boundary fiber sum of $\C$-marked $\RLF s$ is naturally $\C$-marked.

\end{r.}

\begin{p.}\label{iyitanimli}
If $g>1$, then the strong  boundary fiber sum, $X'\natural_{\Sigma_g} X\to
D^2$, of directed $\C$-marked genus-$g$ real Lefschetz fibrations
is well-defined up to isomorphism of $\C$-marked $\RLF s$.
\end{p.}

\noindent {\it Proof:} 
Note that the boundary fiber sum does not affect the fibrations
outside a small neighborhood of the interval where the gluing is made. 
Let us choose a neighborhood $N$ which is real and far from the critical set. Obviously, the real structures on the fibers over the real points of $N$ are isotopic. 
Therefore, each fiber sum defines a path in the space of real structures on $\Sg$, and the difference of two strong boundary fiber sums gives
a loop in this space. Thus, the result follows from the contractibility (shown in Proposition~\ref{fundgrup}) of this loop in the case of $g>1$.  
\hfill  $\Box$ \\


\section{Real Lefschetz chains associated to $\C$-marked real Lefschetz fibrations}\label{gzincirbolum}

\begin{figure}[ht]
   \begin{center}
     \includegraphics[scale=0.29,trim=0 0 180 0]{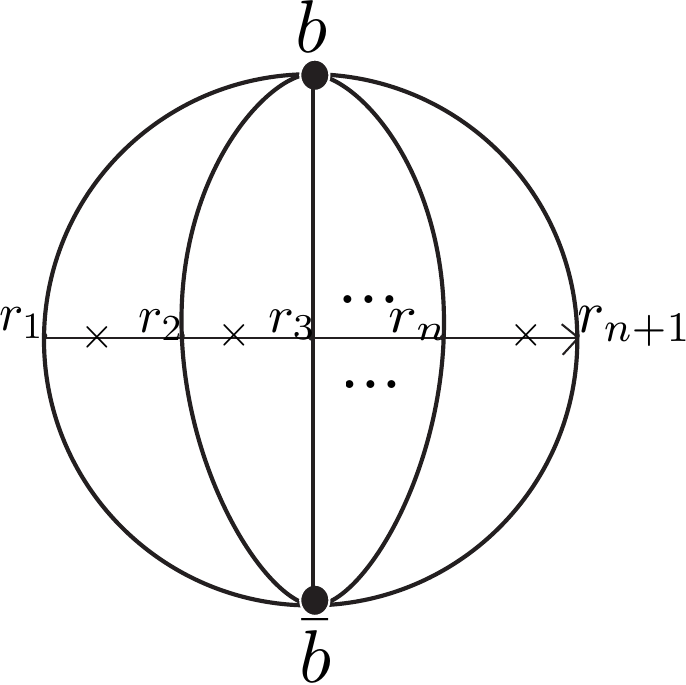}\hspace{1cm}
      \caption{Slicing $D^2$ into small discs having one critical value.}
       \label{karpuz}
       \end{center}
\end{figure}

Let us now consider a directed $\C$-marked real Lefschetz fibration  $\pi:X\to D^2$ with only real  critical values. We slice $D^2$ into smaller discs, $D_1, D_{2},\ldots, D_{n}$  (ordered with respect to the orientation of the real part of $(D^2,\cj)$) such that each $D_i$ contains only one critical value  and the base point $b$ (which is chosen to be the ``north pole", see Figure~\ref{karpuz}).
Let $r_1, r_2,\ldots, r_n, r_{n+1}$ be the real points of $\cup_{i=1}^n \partial D_i$ and let $c_i$ be the  real structure on $\Sg$  pulled back from the inherited real structure of $F_{r_i}$. 

As asserted in Remark~\ref{fcc},  for each fibration over $D_{i}$ we have $c_{i+1}\circ c_{i}=t_{a_i}$ where $a_i$ denotes the corresponding vanishing cycle.  Moreover, as shown in Proposition~\ref{reelkodC}, each $\C$-marked real Lefschetz fibration over $D_i$ is determined by the isotopy class $[c_i, a_i]$ of a real code. Therefore, the fibration $\pi:X\to D^2$  yields a sequence of real codes $[c_i, a_i]$ satisfying $c_{i+1}\circ c_{i}=t_{a_i}$. Obviously this sequence is an invariant of $\pi$. 

\begin{d.}
A sequence $[c_1, a_1], [c_2, a_2],...,[c_n, a_n]$  is called a {\em real Lefschetz chain}  (of isotopy classes of real codes) if
we have $c_{i+1}\circ c_{i}=t_{a_i}$ for all $i=1,...,n$. 
\end{d.}

\begin{t.}\label{gucluzincir} If $g>1$, then there is a one-to-one correspondence between
the real Lefschetz chains $[c_1, a_1], [c_2, a_2],...., [c_n,
a_n]$ and the isomorphism classes of directed
$\C$-marked genus-$g$ real Lefschetz fibrations over $D^2$ with
only real critical values.
\end{t.}

\noindent {\it Proof:}  Necessity is clear. 
As for the converse, we consider the unique class (assured by  Proposition~\ref{reelkodC})  of  directed $\C$-marked elementary real Lefschetz fibration associated to each real code $[c_{i}, a_{i}]$. We then glue these elementary fibrations (from left to right respecting the order determined by the chain) using the strong boundary fiber sum. The result, thus, follows from Proposition~\ref{iyitanimli}.
\hfill  $\Box$ \\

If the total monodromy  of the fibration $\pi: X\to D^2$ is
the identity, then we can consider the extension of $\pi$ to a fibration
$\hat{\pi}:\hat{X}\to S^2$. Two such extensions,
$\hat{\pi}:\hat{X}\to S^2$ and $\check{\pi}:\check{X}\to S^2$,
are considered {\em isomorphic} if there is an equivariant orientation preserving
diffeomorphism $H: \hat{X}\to \check{X}$ such that
$\hat{\pi}=\check{\pi}\circ H$. 

\begin{p.}\label{kureyecik}  Let  $\pi:X \to D^2$  be a $\C$-marked genus-$g$ real Lefschetz fibration whose total monodromy is the identity. 
If $g>1$, then $\pi$ can be extended uniquely up to isomorphism to a real Lefschetz fibration over $S^2$.
\end{p.}

\noindent {\it Proof:} 
Once again, the difference of two  extensions corresponds to a loop in the space of real structures. Hence, the result follows from Proposition~\ref{fundgrup}.
\hfill  $\Box$ \\

\begin{c.}\label{kuredegucluzincir}
If $g>1$, then there is a one-to-one correspondence between
the real Lefschetz chains $[c_1, a_1], [c_2, a_2],...., [c_n,
a_n]$ of the isotopy classes of real codes such that $c_{n+1}\circ c_{1}= ( t_{a_{n}}\circ c_{n})\circ c_{1}=id$ and the isomorphism classes of directed
$\C$-marked genus-$g$ real Lefschetz fibrations over $S^2$ with
only real critical values. \hfill  $\Box$ \\
\end{c.}

\begin{r.} 
It is known that the components  of the space of diffeomorphisms of the torus fixing a point is contractible \cite {EE}, so Theorem~\ref{gucluzincir}  can be adapted  to $\C$-marked real elliptic Lefschetz fibration admitting a real section  (a section compatible with the real structures). Details can be found in \cite[ Section~5.4]{N1}.  In the next  section, we prefer to treat the case of elliptic Lefschetz fibrations which may not admit a real section. Moreover, we concentrate on the case of non-marked fibrations. \end{r.}


\section{Boundary fiber sum of non-marked real elliptic Lefschetz fibrations}\label{ftoplama}
To deal with the case of elliptic fibrations, we introduce the  boundary fiber sum for non-marked fibrations.   (Although we concentrate on the case of $g(F)=1$, the definition applies to any genus.) 

\begin{d.}
Let  $\pi':X'\to D^2 $ and $\pi:X \to D^2 $ be two directed
non-marked $\RLF s$. We consider the real fibers, $F'_+$ and $F_-$ of
$\pi'$ and $\pi$ over the real points $r'_+$ and $r_-$,
respectively. Let us assume that the real structure $c'_+:F'_+\to
F'_+$ is conjugate to $c_-:F_-\to F_-$. Namely, there is
an orientation preserving equivariant diffeomorphism
$\phi:F'_{+}\to F_{-}$. Then,  the {\em  boundary fiber sum} of $X'\natural_{F, \phi} X\to D^2$ is obtained
 by identifying the fibers $F'_{+}$ and  $F_{-}$ via
$\phi$.
\end{d.}

\begin{figure}[ht]
   \begin{center}
  \includegraphics[scale=0.30,trim=0 0 80 0]{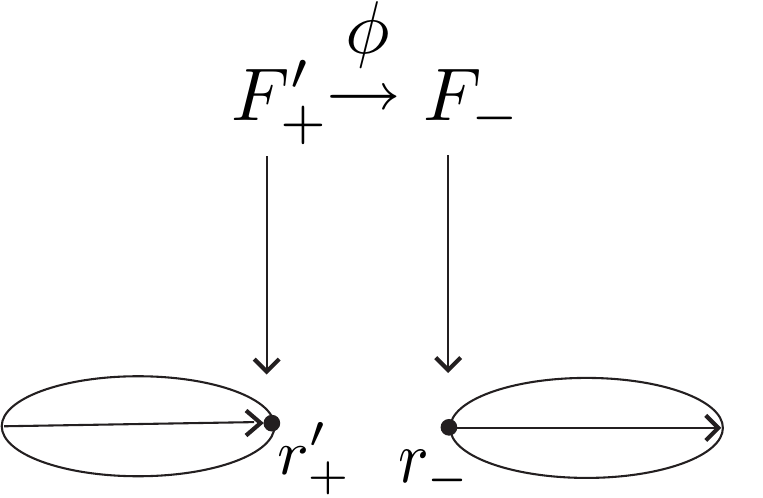}\hspace{1cm}
       \end{center}
 \end{figure}

The boundary fiber sum does depend on the choice of $\phi$ in such a way that
the two boundary fiber sums defined by the equivariant diffeomorphisms
$\phi, \psi : F'_{+}\to F_{-}$ are isomorphic, if
 $\psi \circ \phi^{-1}:F_{-}\to F_{-}$
can be extended to an equivariant diffeomorphism of $X\to D^2$ (or similarly
if $\phi^{-1} \circ \psi:F'_{_+}\to F'_{_+}$ can be extended to an
equivariant diffeomorphism of $X'\to D^2$). 
The necessary and sufficient condition for  $\psi \circ \phi^{-1}: F_{-}\to F_{-}$ to extend to an equivariant diffeomorphism of the fibration $X\to D^2$, is that  $ \psi \circ \phi^{-1}$ takes the unique vanishing cycle $a$ of $X\to D^2$ to a curve equivariantly isotopic to $a$.

Note that  if  $c(a)=a$, then $c$ induces an action on $a$. Such an action can be the identity, a reflection or an antipodal involution.
It is not hard to show that if $c: \Sigma_{1} \to \Sigma_{1}$  has one real component, then
$\Sigma_{1}$ contains a unique $c$-equivariant isotopy class of non-contractible curves on which $c$ acts as a reflection, a unique class of
 curves where the action of $c$ is an antipodal involution, and a unique real curve;
if $c$ has 2 real components, then $\Sigma_{1}$ contains 
no $c$-equivariant isotopy class of curves on which $c$ acts as an antipodal involution, a unique class of
curves on which $c$ acts as a reflection, and two classes of real curves (in which case, we call a pair of representatives of different classes {\em $c$-twin curves}); if $c$ has no real components, then there exist two $c$-equivariant isotopy classes where $c$ acts as an antipodal involution (as above, a pair of representatives of different classes are called {\em $c$-twin curves}) and no classes of other types.
The boundary fiber sum is, therefore,  well-defined unless the real structure $c$ has no real component or  $c$ has two real components one of which is the vanishing cycle $a$.  

Recall that the rotation
$R_{\frac{1}{2}}$ (introduced in the proof of Proposition~\ref{fundgrup}) switches
the $c$-twin curves. Hence, $c$-twin curves can be carried to each other via equivariant diffeomorphisms
although they are not equivariantly isotopic, so in the case of existence of $c$-twin curves, there is an ambiguity in
the definition of the  boundary fiber sum $X'\natural X\to D^2$ (it can be
defined in two ways). To resolve the ambiguity, we should
specify how we identify the $c'_{+}$-twin curves on the fiber $F_+'$ in
$X'$ with the  $c_{-}$-twin curves on the fiber $F_-$ in $X$.
In a certain case, namely if the real structure $c'_{+}$ has two real components and acts on the vanishing cycle $a'$ as a reflection, the problem of switching $c$-twin curves can be eliminated via the transformation  introduced below.

Let $\pi: X\to D^2$ be an elementary directed real elliptic Lefschetz fibration such that  the real structure $c_{+}:F_{+}\to F_{+}$  acts on the vanishing cycle as a reflection. As a result, one of $c_{\pm}:F_{\pm}\to F_{\pm}$  has 1 real component while the other has 2 real components.  Without loss of generality, we can assume that the real structure  $c_{-}$ has 1 real component.  Our aim is to construct a transformation,  $T_{sing}$, of $X$ that does not change the isomorphism class of the fibration $\pi:X\to D^2$ and that is identity over $S_{-}\subset \partial D^2$ and  interchanges the real components of $F_{+}$. To construct $T_{sing}$, we consider the following well known model for elementary elliptic fibrations.

\begin{figure}[ht]
  \begin{center}
     \includegraphics[scale=0.25,trim=0 0 0 0]{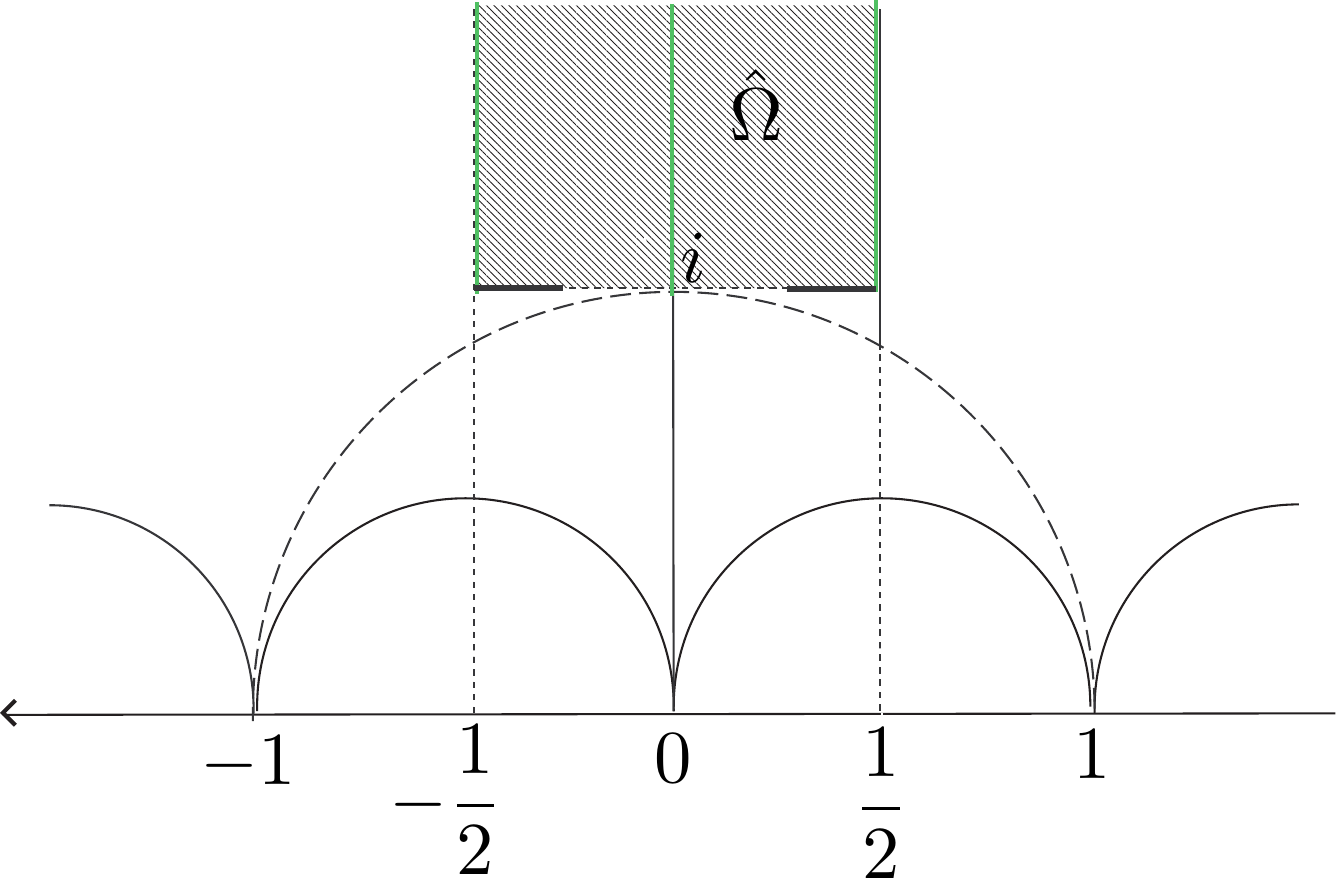}\hspace{0.7cm}
          \includegraphics[scale=0.20,trim=0 -80 00 0]{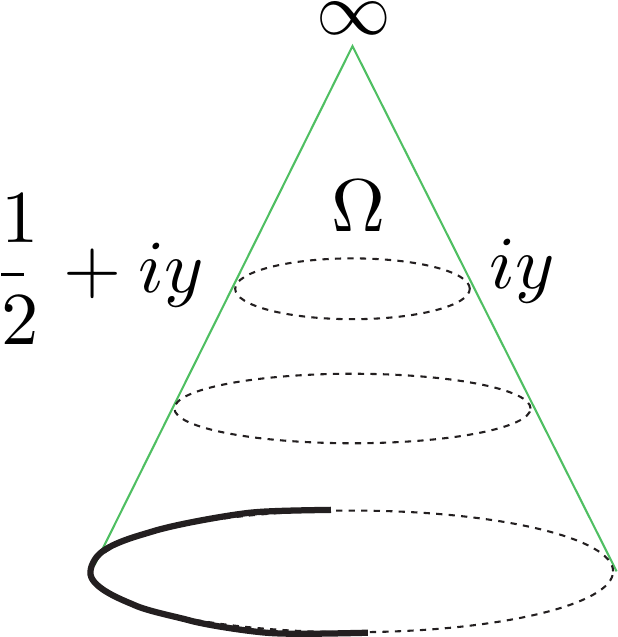}
\caption{Moduli space of prescribed $\RELF s$}
\label{Omega}
      \end{center}
\end{figure}

Let $\hat{\Omega}=\{z| \left|Re(z)\right|\leq \frac{1}{2}, Im(z)\geq 1\}\cup\infty$, (the subset bounded by $Im(z)\geq1$ of the one point compactification of the standard fundamental domain $\{z| \left|Re(z)\right|\leq \frac{1}{2}, \left|z\right|\geq 1\}$ of the modular action on $\C$, see Figure~\ref{Omega}).

We consider the real structure $c_{\hat{\Omega}}:{\hat{\Omega}}\to {\hat{\Omega}}$ such that $c_{\hat{\Omega}}(\omega)=\overline{-\omega}$. Let $\Omega$ denote the quotient $\hat{\Omega}\diagup_{\tiny{\frac{1}{2}+iy\sim -\frac{1}{2}+iy}}$. The real structure $c_{\hat{\Omega}}$ induces a real structure on $\Omega$. Note that $\Omega$ is a topological real disc and can be identified with $D^2$ so that the real part of $D^2$ corresponds to the union of the half-lines $iy$ and $\frac{1}{2}+iy$ where $y\geq1$. For any $\omega\in \Omega$, the fiber over $\omega$ is given by $F_\omega= \C/(\Z+\omega\Z)$, where the fiber $F_{\infty}$ has the required nodal type singularity.

Let $\pi_\Omega: X_\Omega\to \Omega$ denote the fibration such that $\pi_\Omega^{-1}(\omega)=F_\omega=\C/(\Z+\omega \Z)$,  $\forall \omega\in \Omega$. Then, we consider the translation $T'_{\Omega}$ defined by

$$\begin{array}{llll}
T'_{ \Omega}: &X_{\Omega}&\to& X_{\Omega}\\
& (z)_{\Z+\omega \Z}\in F_\omega &\to&(z+\tau(w))_{\Z+\omega \Z} \in F_\omega
\end{array}$$
 
where $(.)_{\Z+\omega \Z}$ denotes the equivalence class in $\C/(\Z+\omega\Z)$. 

\begin{figure}[ht]
\begin{center}
\includegraphics[scale=0.30,trim=0 0 0 0]{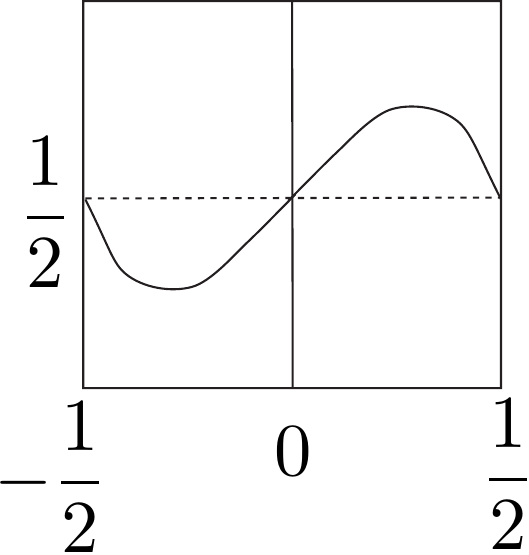}
\caption{The graph of $f$.}
\label{f}
 \end{center}
\end{figure}

We consider  $\tau:\Omega\to \Omega$ such that 
$$\tau(\omega)=-\frac{1}{2}+(\frac{1}{2}-f(Re(\omega))+i)exp(-Im( \omega)+1)$$ where $f:\R/\Z\to \R/\Z$  is a smooth mapping whose graph is as shown in Figure~\ref{f}  and which satisfies the following properties:\\[2mm]

\noindent $\bullet$ $f(0)=\frac{1}{2}$ (modulo $\Z$),\\
$\bullet$ $f(1-x)=1-f(x)$, ($\Rightarrow f(\frac{1}{2})=\frac{1}{2}$) (modulo $\Z$),\\
$\bullet$$f$ is linear on $[\frac{1}{4}, \frac{3}{4}]$ (modulo $\Z$).\\

Note that $\tau$ has the following properties. (Equations are
considered modulo the relation $-\frac{1}{2}+iy\sim
\frac{1}{2}+iy$, $y\geq1$.)

\textbullet{} $\tau(\overline{-\omega})=\overline{-\tau(\omega)},$ \\[3mm]
\indent\textbullet{} $\tau(\infty)=\frac{1}{2},$ \\[3mm]
\indent \textbullet{} $\tau(\frac{1}{2}+iy)=-\frac{1}{2}+iexp(-y+1)=\frac{1}{2}+iexp(-y+1),$\\[2mm]
in particular, if $y=1$, then $\tau(\frac{1}{2}+i)= \frac{1}{2}+i$,\\[2mm]
\indent \textbullet{}  $\tau(iy)=-\frac{1}{2}+iexp(-y+1)=\frac{1}{2}+iexp(-y+1),$\\[2mm]
in particular, if $y=1$, then $\tau(i)=\frac{1}{2}+i$.

\begin{figure}[ht]
  \begin{center}
     \includegraphics[scale=0.39,trim=0 0 170 0]{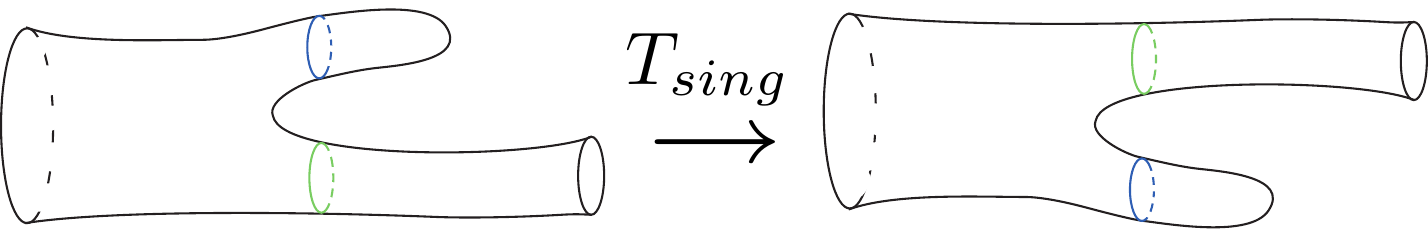}\hspace{1cm}
\caption{\small{The action of $T_{sing}$ on the real part.}}
\label{Tnsing}
      \end{center}
\end{figure}

Let $T_{sing}:X\to X$ denote the transformation induced from $T'_{sing}:X_{\Omega}\to X_{\Omega}$. By definition $T_{sing}$ is equivariant and the identity over $S_{-}\subset \partial D^2$,  and its restriction to $F_{+}$ is the rotation $R_{\frac{1}{2}}$. 
(Figure~\ref{Tnsing} shows the action of $T_{sing}$ on the real part. )

\begin {l.} \label{lifleridondur}
Let $\pi':X'\to D^2$ and $\pi:X \to D^2$ be two non-marked elementary $\RELF s$ such that both $c'_+$ and $c_-$ have 2 real components. We assume that the vanishing cycle $a$ of $\pi$ is real with respect to $c_-$. Then, the  boundary fiber sum $X'\natural_F X \to D^2$ is well-defined if  $c'_+$ acts on the vanishing cycle $a'$  as a reflection.\end{l.}

\noindent {\it Proof:} 
The boundary fiber sums $X'\natural_{F, \phi} X\to D^2\,\, \textrm{and}\,\, X'\natural_{F, \psi} X\to D^2$ are not isomorphic if $\phi\circ \psi^{-1}(a)$ and $a$ are $c$-twin curves. 
However,  in the case when  $c'_+$ acts on the vanishing cycle $a'$  as a reflection, we can apply $T_{sing}$ to $X'$ so that $T_{sing}(F'_+)$ differs from the fiber $F'_+$ by the rotation $R_{\frac{1}{2}}$. Therefore, $X' \natural_{F,\phi} X\to D^2$ is isomorphic to  $T_{sing}(X')\natural_{F,\phi \circ R_\frac{1}{2}} X\to D^2$ which is isomorphic to $X'\natural_{F, \psi} X\to D^2$. 
\hfill  $\Box$ \\


\section{Real Lefschetz chains associated to non-marked real elliptic Lefschetz fibrations}\label{zzincirbolum}

\begin{figure}[ht]
   \begin{center}
      \includegraphics[scale=0.29,trim=0 0 170 0]{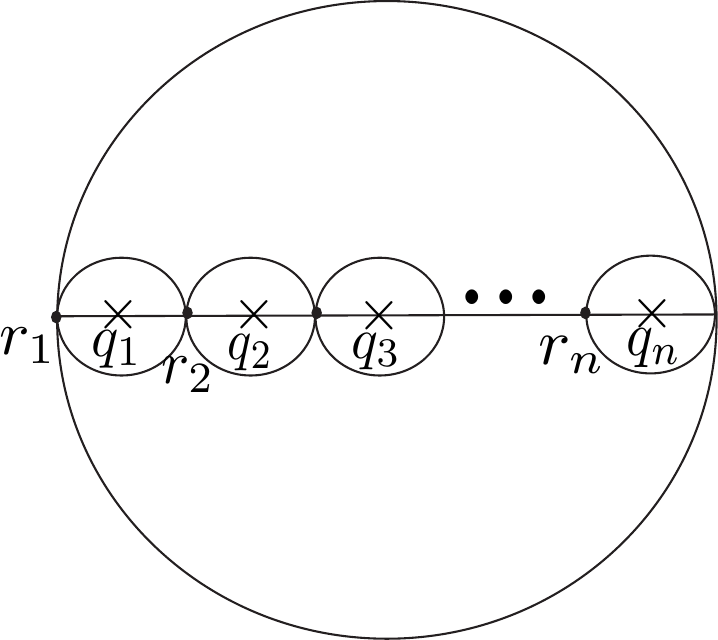}\hspace{1cm}
\caption{Subdividing $D^2$ into smaller discs.}
\label{zayif}
 \end{center}
 \end{figure}

We now consider a non-marked directed real elliptic Lefschetz fibrations $\pi:X\to D^2$ with
only real critical values, $q_1< q_2<...<q_n$. Around each critical value $q_i$ we choose a small real disc $D_i$ such that
$D_i\cap \{q_1, q_2,...,q_n\}=\{q_i\}$ and $D_i \cap D_{i+1}=\{r_{i+1}\} \subset [q_i, q_{i+1}]$, see Figure~\ref{zayif}. 
Let $c_{i}$ be the real structures on the fibers $F_{r_{i}}$, $1\leq i\leq n$ (where $r_{1}$ is the left real point of $\partial D^2$) and $a_{i}$ be the corresponding vanishing cycle.  

By  Proposition~\ref{reelkod}, each directed (non-marked) fibration over $D_i$ is classified by the conjugacy class $\{c_i, a_i\}$ of the
real code. Thus, we can encode the fibration $\pi:X\to D^2$ by the sequence 
$\{c_1, a_1\}, \{c_2, a_2\},..., \{c_n, a_n\}$ of conjugacy classes of real codes such that $ t_{a_i} \circ c_{i}$ is
conjugate to $c_{i+1}$ for all $1\leq i\leq n$.  We call this
sequence the {\em real Lefschetz chain} (of conjugacy classes of real codes).

\begin{figure}[ht]
   \begin{center}
      \includegraphics[scale=0.28,trim=0 0 -10 0]{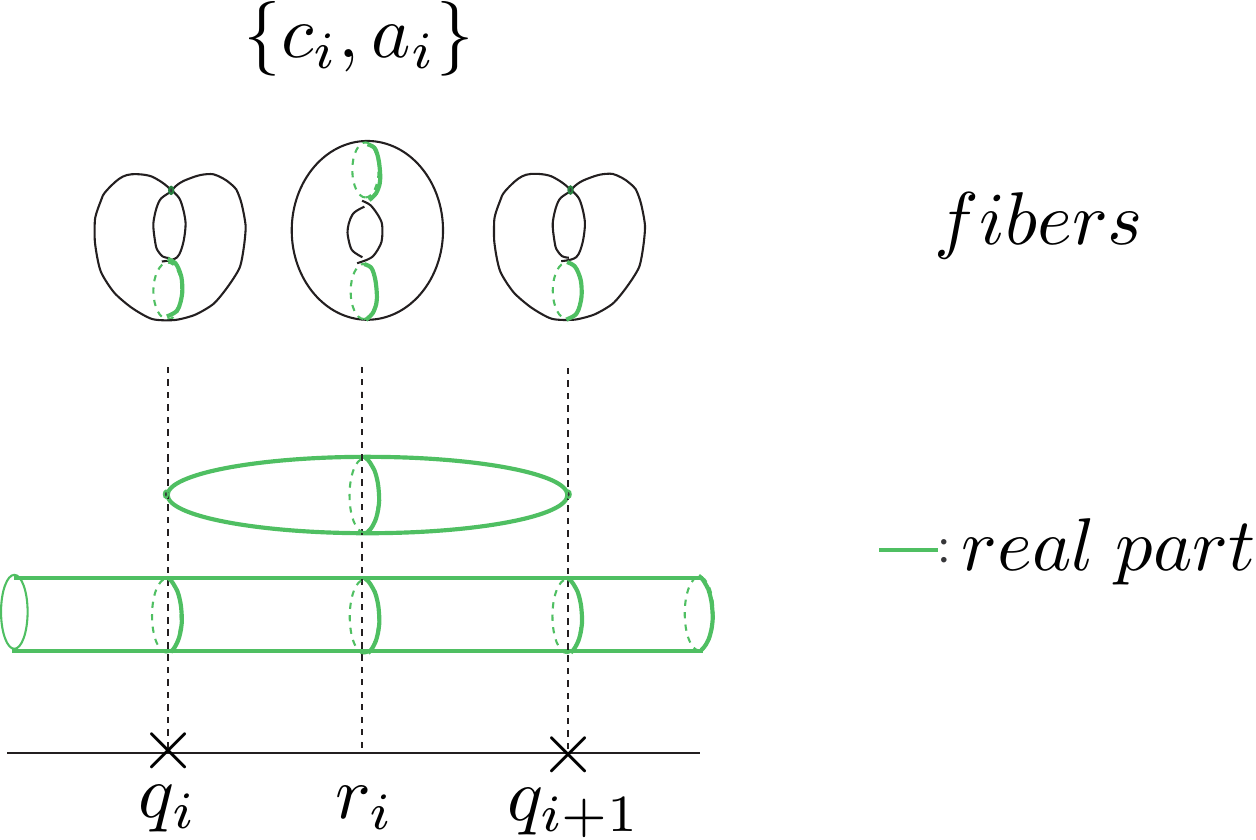}\hspace{.2cm}
   \includegraphics[scale=0.28,trim=0 0 -10 0]{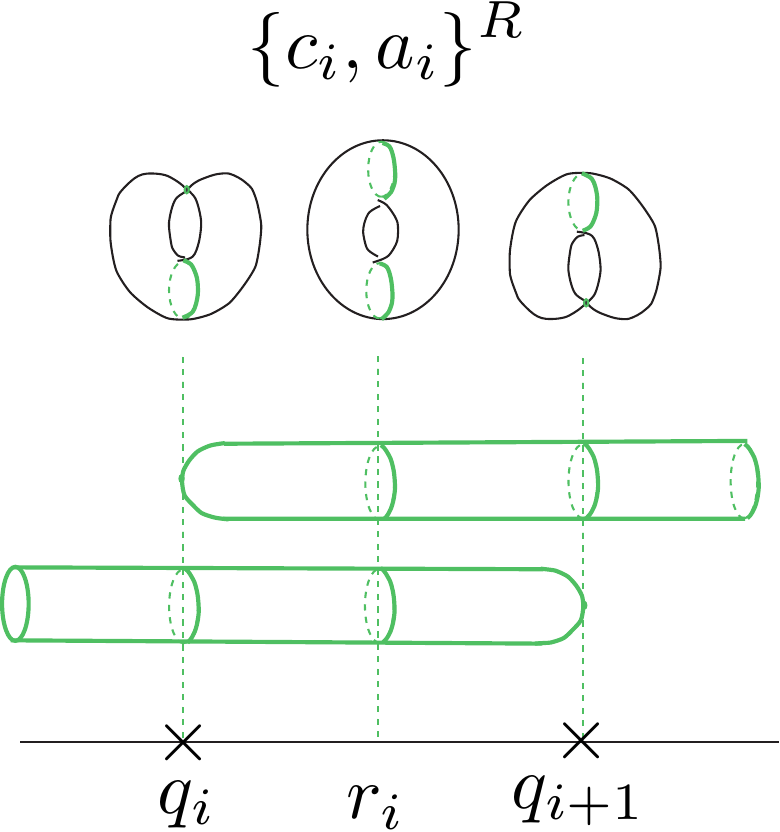}\vspace{.7cm}
  \includegraphics[scale=0.28,trim=0 0 10 0]{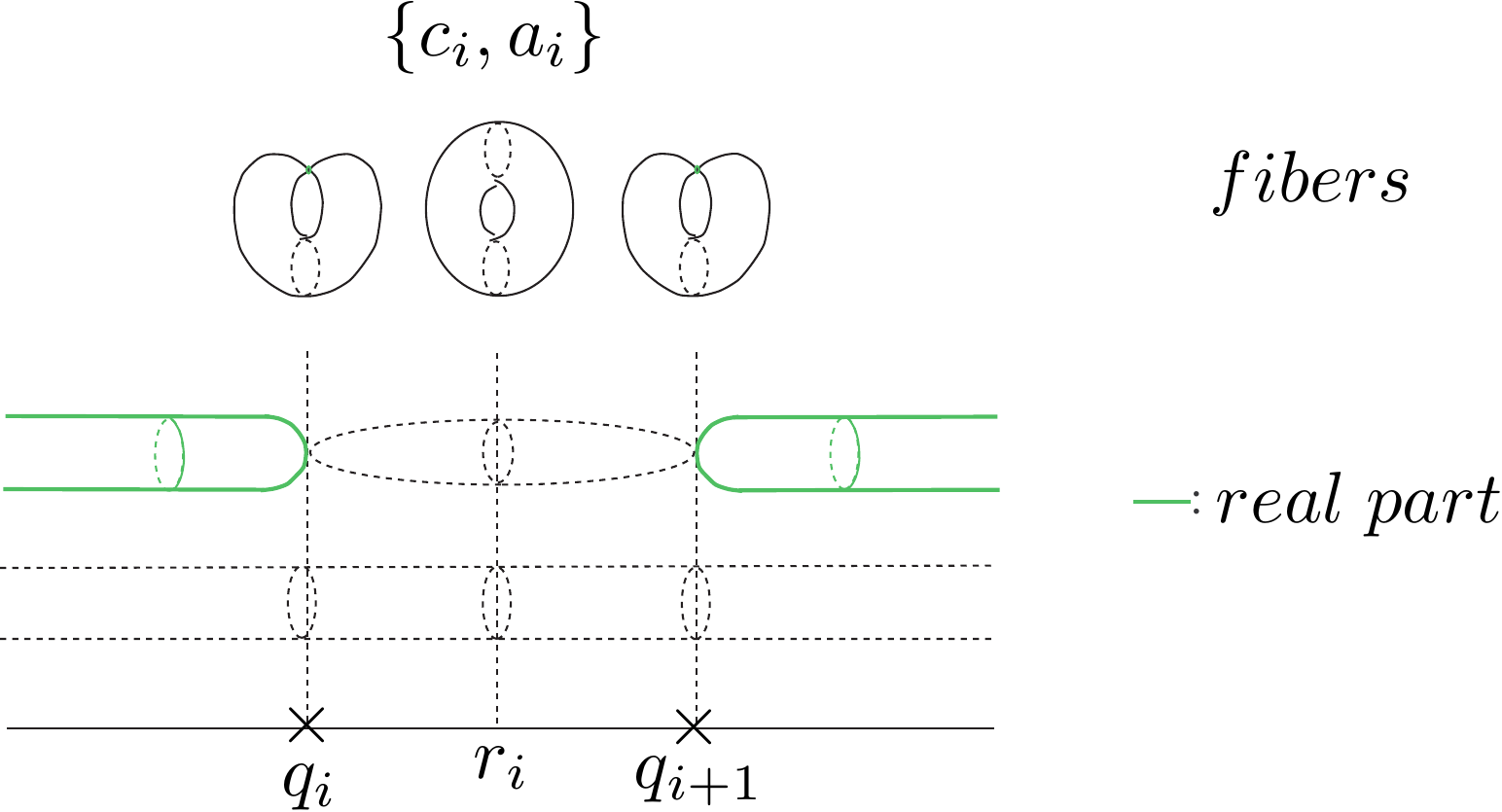}\hspace{.4cm}
  \includegraphics[scale=0.28,trim=0 0 10 0]{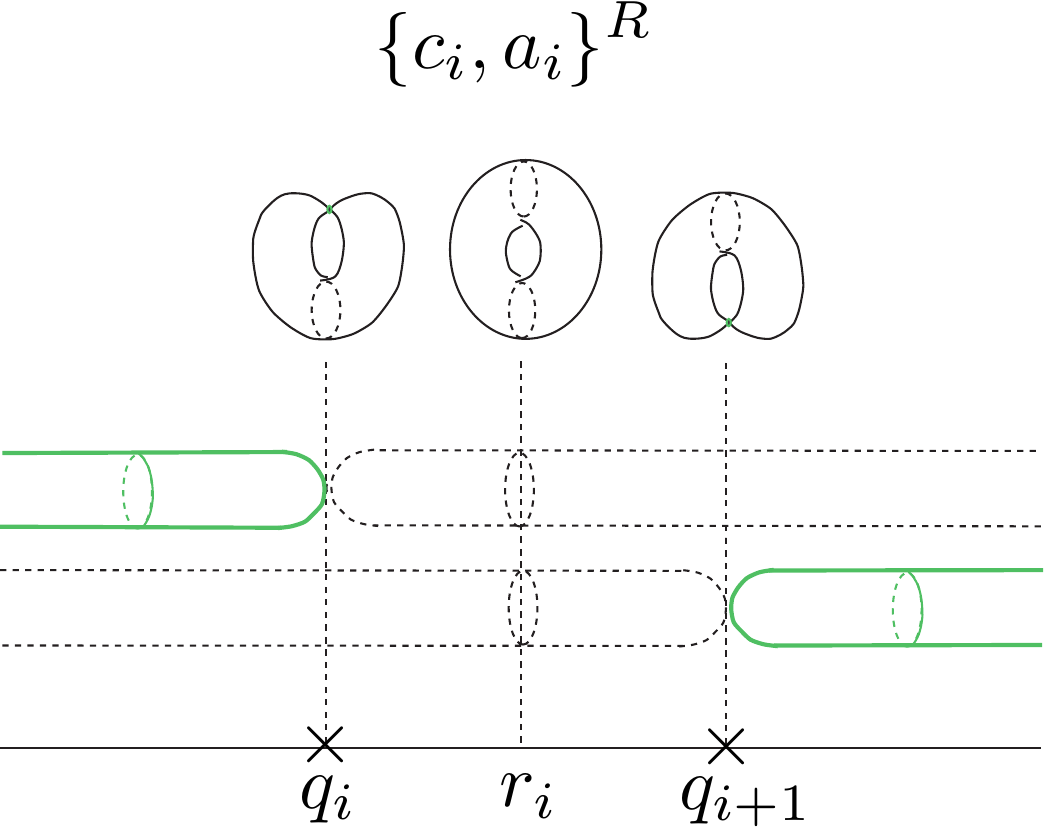}
\caption{ Real parts of the fibrations associated to $\{c_i, a_i\}$ and $\{c_i, a_i\}^R$.}
\label{cacar}
       \end{center}
 \end{figure}

Clearly, real Lefschetz chains are invariants of directed non-marked $\RELF s$ over disc that have only real critical values, but they are not sufficient for classifying such fibrations. Additional information is needed, if for some $i$ the real structure  $c_i$ has 2 real components and vanishing cycles corresponding to the critical values $q_i$ and $q_{i+1}$ are real, respectively, or if $c_{i}$ has no real component.  Indeed, in these cases the vanishing cycles corresponding to the critical values $q_i$ and $q_{i+1}$ can be the same curve,  or  they can be $c_i$-twin curves. If they are $c_{i}$-twin curves, then we mark  $\{c_i, a_i\}^R$ 
the corresponding real code $\{c_i, a_i\}$ by  adding $R$   (here $R$ refers to the rotation $R_{\frac{1}{2}}$ which interchanges $c$-twin curves). 
The real Lefschetz chain we obtain is called the {\em decorated real Lefschetz chain}. 
Figure~\ref{cacar} shows all possible configurations of the real locus associated to $\{c_i, a_i\}$ and $\{c_i, a_i\}^R$.

\begin{t.}\label{zayifzincir}
There exists a one-to-one correspondence between the decorated real Lefschetz chains and the isomorphism classes of directed non-marked real elliptic Lefschetz fibrations over $D^2$ that have only real critical values.
\end{t.}

\begin{figure}[ht]
   \begin{center}
  \includegraphics[scale=0.28,trim=0 0 10 0]{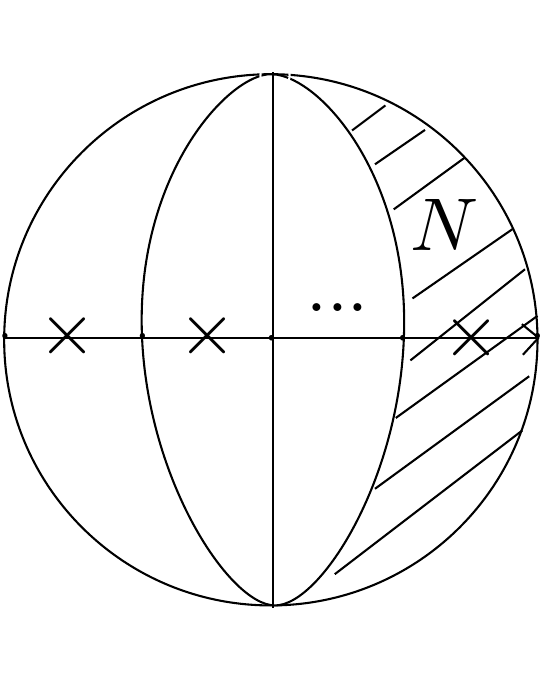}
\caption{Neighborhood over which $T_{sing}$ is applied.}
\label{N}
       \end{center}
 \end{figure}

\noindent {\it Proof:} 
Necessity is clear. As for the converse, we consider the unique class of directed non-marked elementary $\RELF$ (assured by  Proposition~\ref{reelkod}) associated to each real code $\{c_i, a_i\}$.
Then, we construct the required fibration by gluing elementary fibrations (from left to right) using the boundary fiber sum. 
As is discussed above, the boundary fiber sum is uniquely defined in the case when the real structure on the fiber where the sum is performed has 1 real component or when it has 2 real components and acts on the vanishing cycle of the elementary fibration glued to right as a reflection. In the case when the real structure has 2 real components and acts on the last vanishing cycle of the already constructed fibration  $\pi':X'\to D^2$ as  a reflection,  the two possible  boundary fiber sums are isomorphic by Lemma~\ref{lifleridondur} since in this case we can apply $T_{sing}$ to $X'$ (by considering $T_{sing}$ on a neighborhood $N$ of the last critical value, as shown in Figure~\ref{N}, and extending it to $X'$ as the identity outside of $\pi'^{-1}(N)$). In all the other cases, the boundary fiber sum is defined uniquely by the decoration.
\hfill  $\Box$ \\

If $c_1$ is conjugate to $c_{n+1}$, then we can consider an extension of $\pi: X\to D^2$ to a fibration over $S^2$.
As before, in the case when $c_{n+1}$ has no real components  or it has 2 real components and both $a_1$ and $a_n$  are real, a decoration at infinity will be needed.

\begin{p.}\label{kureyeyay} 
Let $\pi: X\to D^2$ be a real elliptic Lefschetz fibration associated to a decorated real Lefschetz chain.  We assume that the real structures $c_{1}$ and $c_{n+1}$ on the fibers over left and respectively right real point of $\partial D^2$ are conjugate. 
If  $c_{n+1}$ (and thus $c_{1}$) has 1 real component or if $c_{n+1}$ (and thus $c_{1}$) has 2 real components and either $c_{n+1}$  acts on the vanishing cycle $a_{n}$ as a reflection, or  $c_{1}$ acts on the vanishing cycle  $a_{1}$ 	as a reflection, then $\pi$ extends uniquely  to a fibration over $S^2$. 
Otherwise, there are two extensions distinguished by the decoration at infinity.
\end{p.}

\noindent {\it Proof:} 
An extension of $\pi:X\to D^2$ to a fibration over $S^2$ defines a trivialization, $\phi:\Sa\times S^1 \to \pi^{-1}(\partial D^2)$ over the boundary $\partial D^2$. Two trivializations $\phi, \phi'$ correspond to isomorphic real fibrations if $\phi^{-1} \circ \phi': \Sa\times S^1\to \Sa\times S^1$ can be extended to an equivariant diffeomorphism of $\Sa\times D^2$ with respect to the real structure $(c_{n+1}, \cj):\Sa\times D^2 \to \Sa\times D^2$.
Let ${\Phi}_t=(\phi^{-1} \circ \phi')_t:\Sa\to \Sa$, $t\in S^1$. Since there is no fixed marking, up to change of marking we assume that ${\Phi}_t\in \Do(\Sa)$.

The real structure splits the boundary into two symmetric pieces, so instead of considering an equivariant map over the entire boundary we consider a diffeomorphism  over one the symmetric pieces. Let $\Phi_t, t \in [0,1]$ denote the family of such diffeomorphisms. The family, $\Phi_t$, $t\in [0,1]$ defines a path in $\Do(\Sa)$ whose end points lie in the group $\Dcno(\Sa)$; therefore, $\Phi_t$ defines a relative loop in $\pi_1(\Do(\Sa),\Dcno(\Sa))$, and we are interested in the contractibility of this relative loop. 

We consider the following exact sequence of the pair $(\Do(\Sa),\Dcno(\Sa))$

$...\to\pi_1(\Dcno)\to \pi_1(\Do)\stackrel{f}{\rightarrow} \pi_1(\Do, \Dcno)\stackrel{g}{\rightarrow} \pi_0(\Dcno)\stackrel{h}{\rightarrow}$\\
$ \pi_0(\Do)\to \pi_0(\Do, \Dcno)\to 0.$

In the case when $c_{n+1}$ has one real component, $\Dcno(\Sa)$ is connected, so the map $h$ is injective, so $f$ is surjective. Therefore, elements of the group $\pi_1(\Do(\Sa), \Dcno (\Sa), id)$ can be seen in $\pi_1(\Do(\Sa), id)$. 

In all the other cases, $\Dcno(\Sa)$ has two components.  We mark one of the components to get  the map $h$ injection, when restricted to the marked component. Thus, $g$ becomes the zero map, and so $f$ is surjective over the marked component of $\Dcno(\Sa)$. Note that decoration of real Lefschetz chains distinguishes one of the component of $\Dcno(\Sa)$; hence, marking one component or other give the two extensions distinguished by the decoration.  The  distinctive feature of the case when $c_{n+1}$ has 2 real components and acts  $a_n$ as a reflection (or  $c_{1}$ acts on $a_{1}$ as a reflection) is that  the transformation $T_{sing}$ changes one marking to other, so the marking is not essential. 

The proposition, thus,  follows from Lemma~\ref{bazN}  in which we show that any relative loop 
can be made contractible by means of some transformations $T$ of the fibration $\pi:X\to D^2$.
\hfill  $\Box$ \\

Let us first define the transformation $T$ of real elliptic Lefschetz fibrations  over $D^2$ that is defined over a regular slice $N$ of $D^2$.

\begin{figure}[ht]
  \begin{center}
     \includegraphics[scale=0.28,trim=0 0 170 0]{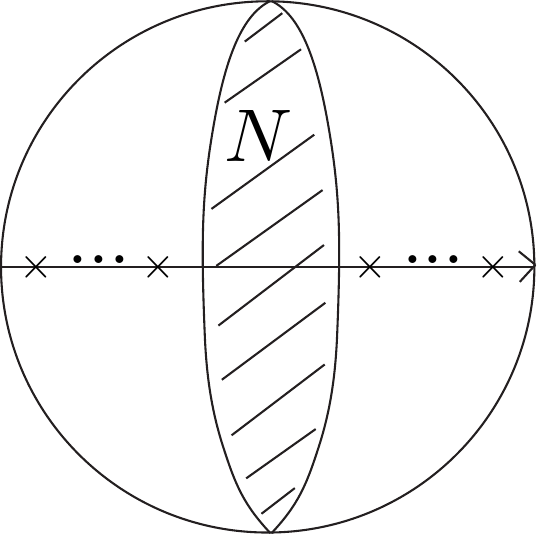}\hspace{1cm}
\caption{Neighborhood over which $T$ is applied.}
\label{3N}
      \end{center}
\end{figure}

Let $\pi:X\to D^2 $ be a directed $\RELF$. We consider a real slice $N$ of $D^2$ which contains no critical value, see Figure~\ref{3N}.

Let $\xi:I\times I \to N$, $I=[0,1]$ be an orientation preserving diffeomorphism such that first interval correspond to the real direction on $N$.
The fibration over $N$ has no singular fiber; hence, it is trivializable. Let us consider a trivialization $\Xi: \Sa\times I \times I \to \pi^{-1}(N)$ such that the following diagram commutes
$$\xymatrix{\Sa\times I\times I \ar[r]^{\Xi} \ar[d] & \pi^{-1}(N) \ar[d]^{\pi} \\
          I\times I \ar[r]^{\xi} & N.}$$

Since $N$ has no critical value, the isotopy type of the real structure on the fibers over the real part of $N$ remains fixed. 
If the real structure $c$ has 2 real components, then we consider the model $\varrho:\C/\Z^2\to \Sa$
and set  $\bar{\varrho}=(\varrho, id):  \C/\Z^2 \times I \times I\to \Sa \times I\times I$   to define $T$ as follows 

$$\begin{array}{llll}
T':&\C/\Z^2 \times I\times I&\to& \C/\Z^2 \times I\times I\\
& ((x+iy)_{\Z^2}, t, s)&\to&((x+ t+iy)_{\Z^2}, t, s).
\end{array}$$

Then, we set  $T= \Xi\circ (\bar{\varrho}\circ T' \circ \bar{\varrho}^{-1}) \circ \Xi^{-1}$ on $\pi^{-1}(N)$.
Since $T$ is the identity at $t=0, 1$, we can extend $T$ to $X$ by the identity outside of $\pi^{-1}(N)$. 

If $c$ has 1 real component, then we construct the transformation $T$ using $\varrho: \C/ \Lambda\to \Sa$. Similarly, if $c$ has no real component, then we repeat the same using $\varrho:\R^2/ \Z^2\to \Sa$.

\begin{r.}
1. Since the transformation $T$ is defined by a real rotation, $T$ preserves the isomorphism class of the real Lefschetz fibration.

2. The map $T$ depends only on the isotopy type of $\pi^{-1}(N)$.
\end{r.}

\begin{l.}\label{bazN}
Let $\pi:X\to D^2$ be a real elliptic Lefschetz fibration with only real critical values. We assume that there exists at least one vanishing cycle on which corresponding real structure acts as a reflection.
Then, there exists a generating set for $\pi_1(\Do(\Sa), id)=\Z+\Z$ consisting of transformations $T_{\pm}$ for some non-singular slices $N_{\pm}$.
\end{l.}

\noindent {\it Proof:} 
Let $q_{i}$ be the critical value such that the real structure on a nearby regular real fiber acts on the vanishing cycle as a reflection. 
This assumption assures that the neighboring real fibers have one real component on one side and 2 real components on the  other side of the critical value $q_{i}$. Without loss of generality we can assume that the real structure over a fiber over a  real point which lies on the left of  $q_i$ has 2 real components. (The other case can be treated similarly.)

\begin{figure}[ht]
  \begin{center}
       \includegraphics[scale=0.32,trim=0 0 70 0]{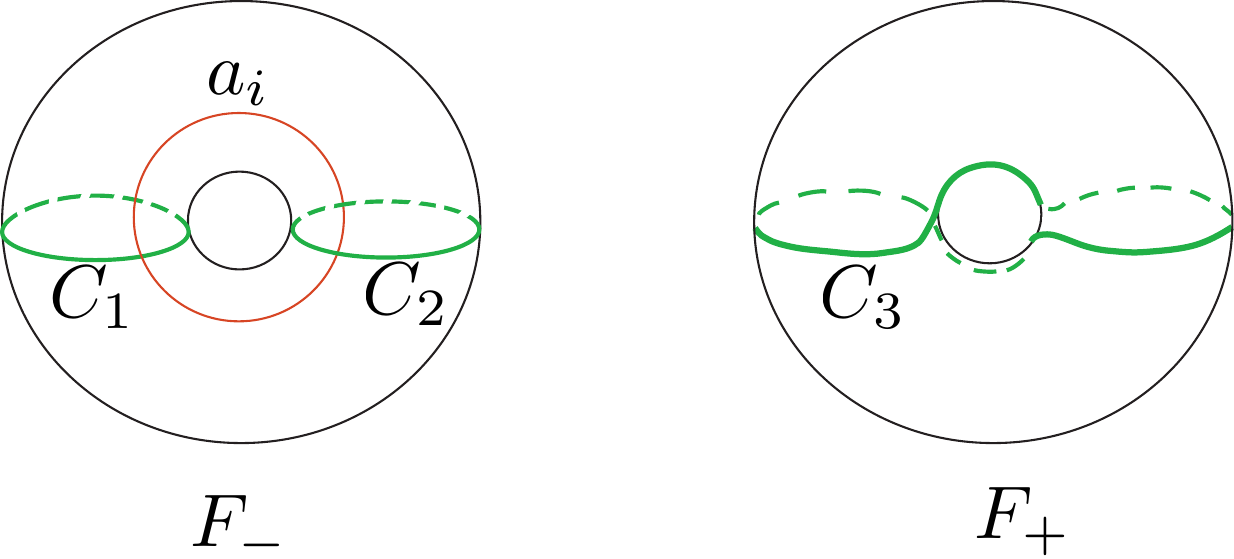}
\caption{Real fibers over the real points neighboring $q_{i}$.}
\label{2aeksib}
      \end{center}
\end{figure}

We choose an auxiliary $\C$-marking $(\{b, \bar{b}\}, \{ \rho:\Sa\to F_b , \bar{\rho}: \Sa \to F_{\bar{b}}\})$ and fix an identification $\varrho: S^1\times S^1\to \Sa$. Since the real structure has 2 real components, we can assumed that the induced real structure on $S^1\times S^1$ is the reflection $(\alpha, \beta)\to (\alpha, -\beta)$. The real part consists of the curves $C_1=(\alpha, 0)$ and $C_2=(\alpha, \pi)$. Moreover, a representative of the vanishing cycle can be chosen as $(0, \beta)$. As $c_+=t_{a_i}\circ c_-$ on $S^1\times S^1$ the real part of $c_+$ is the curve, $C_3$, given  homologically by $2\alpha-\beta$ ( see Figure~\ref{2aeksib}).

\begin{figure}[ht]
  \begin{center}
        \includegraphics[scale=0.27,trim=0 0 70 0]{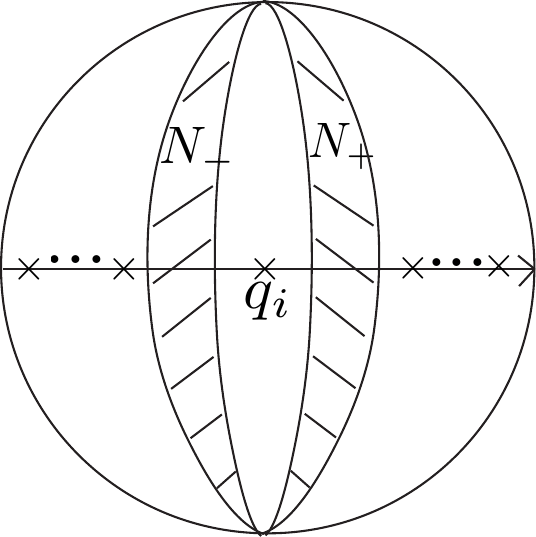}
\caption{Regular slices $N_{\pm}$.}
\label{N+-}
      \end{center}
\end{figure}

We now consider two non-singular real slices $N_-$, $N_+$  of $D^2$  as shown in Figure~\ref{N+-}.  Let us suppose that the real fibers over $N_{-}$  are identified to $F_{-}$  while real fibers over $N_{+}$ are identified to $F_{+}$ (where $F_{\pm}$ are as shown Figure~\ref{2aeksib}). Let $C'_{3}$ and $C'_{1}$ be curves on $F_{b}$ obtained by pulling back $C_3\subset F_{+}$ and  $C_1\subset F_{-}$, respectively. The curves   $C'_{3}$ and $C'_{1}$ intersect at one point, so we can identify $\Sa$ with $C'_1\times C'_3$ so that rotations along $C'_1$ and $ C'_3$ generate the group $\Do(\Sa,id)$.
Hence, $\{T_{+}, T_{-}\}$ generates $\pi_1(\Do(\Sa), id)$.
\hfill  $\Box$ \\

Theorem~\ref{zayifzincir}  applies naturally to directed non-marked $\RELF s$ over $D^2$ which admit a real section in which case real Lefschetz chain does not contain a real code $(c_i, a_i)$ where the real structure has no real component. Besides, in the case when the real structure has 2 real components and the vanishing cycle is real, the decoration is not needed since the existence of a real section determines naturally the gluing.  Moreover, the extension to a fibration over $S^2$ is uniquely defined by  the section. Hence we have the following proposition. 

\begin{p.} \label{kilciklizayif}
Two directed $\RELF s$ over $S^2$ admitting a real section and having the same real Lefschetz chain up to cyclic ordering are isomorphic.
$\Box$
\end{p.}

\begin{r.} Indeed, the proposition holds even for fibrations with a fixed real section. If there are only real critical values, then the real sections are determined in a neighborhood of a real part. Moreover, over the real part one can carry one real section to another using the transformations $T$  and  {\em double} $T_{sing}$. Indeed, the {\em double} $T_{sing}$  is defined for real Lefschetz fibrations with two critical values where the real structure extracted from the real fiber over a real point between the critical values acts on the vanishing cycles as a reflection. The model we use to define the {\em double} $T_{sing}$ is as follows. Consider the disc $D$ with two critical values as the double cover of a disc with one critical value  branched at a regular real point. Let $D_{-}$ and $D_{+}$ be two corresponding copies of the disk on the branched cover. By pulling back the fibration $X_{\Omega}$ over $D$ we obtain a model fibration over $D_{-}\cup D_{+}$. Thus, we can apply $T_{sing}$ at the same time to fibrations over $D_{-}$ and $D_{+}$.    
The possible modifications of the section is shown in the Figure~\ref{Tnsektdegis}.

\begin{figure}[ht]
  \begin{center}
     \includegraphics[scale=0.25,trim=0 0 170 0]{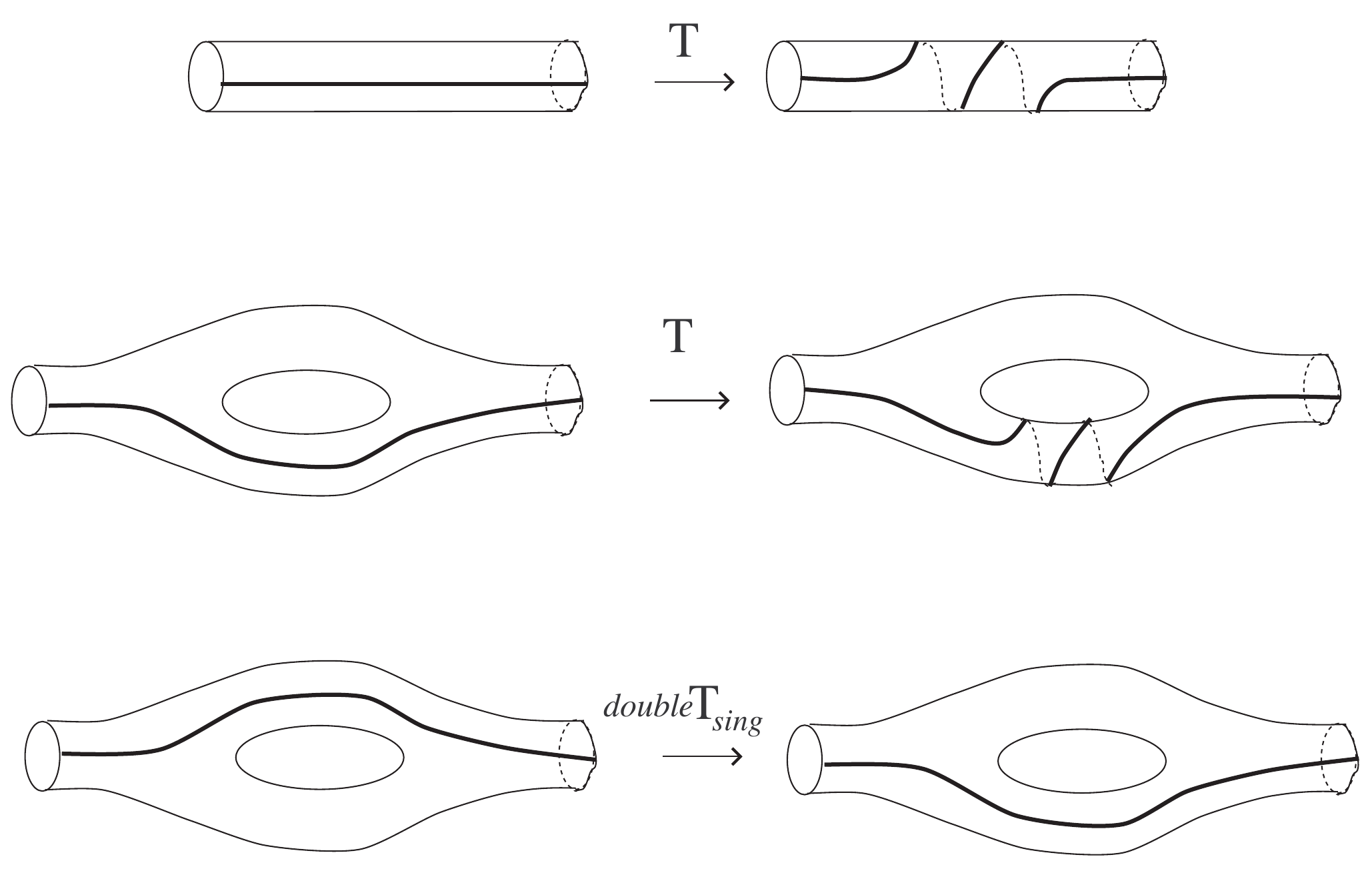}\hspace{1cm}
\caption{\small{Modification of the real section over the real part.}}
\label{Tnsektdegis}
      \end{center}
\end{figure}
\end{r.}


\end{document}